\documentclass[10pt]{article}
\usepackage{amsmath}
\usepackage{amsfonts}
\usepackage{amssymb}
\usepackage{xcolor}

\providecommand{\U}[1]{\protect\rule{.1in}{.1in}}
\begin{document}

\title{An optimal spectral inequality for degenerate operators}
\author{R\'{e}mi Buffe\thanks{
Universit\'e de Lorraine, CNRS, Inria, IECL, F-54000 Nancy, France} , Kim Dang Phung\thanks{%
Institut Denis Poisson, Universit\'{e} d'Orl\'{e}ans, Universit\'{e} de
Tours \& CNRS UMR 7013, B\^{a}timent de Math\'{e}matiques, Rue de Chartres,
BP. 6759, 45067 Orl\'{e}ans, France E-mail address: kim\_dang\_phung@yahoo.fr%
} , Amine Slimani \thanks{Ecole des Mines de Nancy, Université de Lorraine, Campus Artem, CS 14 234, 92 Rue Sergent Blandan, 54042 Nancy}}
\date{26/09/2023 \ }
\maketitle

\bigskip

Abstract .- In this paper we establish a Lebeau-Robbiano spectral inequality
for a degenerate one dimensional elliptic operator. Carleman techniques and
moment method are combined. Application to null controllability on a
measurable set in time for the degenerated heat equation is described.

\bigskip

\bigskip

\section{Introduction and main results}

\bigskip

\bigskip

The purpose of this article is to prove a spectral inequality for a family
of degenerate operators acting on the interval $\left( 0,1\right) $. In
arbitrary dimension, for a second-order symmetric elliptic operator ${%
\mathcal{P}}$ on a bounded domain $\Omega $ with homogeneous Dirichlet or
Neumann boundary conditions, the spectral inequality also called
Lebeau-Robbiano estimate takes the form%
\begin{equation}
\left\Vert u\right\Vert _{L^{2}\left( \Omega \right) }\leq ce^{c\sqrt{%
\lambda }}\left\Vert u\right\Vert _{L^{2}\left( \omega \right) }\text{ ,
\quad }\forall u\in \text{span}\left\{ \Phi _{j};\lambda _{j}\leq \lambda
\right\}  \tag{1.1}  \label{1.1}
\end{equation}%
where $\omega \subset \Omega $ is an open subset and where the functions $%
\Phi _{j}$ form a Hilbert basis of $L^{2}\left( \Omega \right) $ of
eigenfunctions of ${\mathcal{P}}$ associated with the nonnegative
eigenvalues $\lambda _{j}$ , $j\in \mathbb{N}$, counted with their
multiplicities. In other words, the family of spectral projectors associated
to ${\mathcal{P}}$ enjoys an observability inequality on a set $\omega
\subset \Omega $ for low frequencies $\lambda _{j}\leq \lambda $ with a
constant cost as $ce^{c\sqrt{\lambda }}$.

\bigskip

The state of art to prove (\ref{1.1}) is either Carleman inequalities for
elliptic equations (see \cite{LR}, \cite{LZ}, \cite{JL}, \cite{L}, \cite{LRL}%
, \cite{LRLeR1}, \cite{LRLeR2}, \cite{Le}, \cite{LL} and \cite{Q}, \cite{FQZ}%
) or observation estimate at one point in time for parabolic equations (see 
\cite{AEWZ}, \cite{BaP} and \cite{BP}).

\bigskip

One of the key applications of (\ref{1.1}) is either observability for
parabolic systems or controllability for parabolic systems, knowing that
both are equivalent properties by a duality argument (see \cite{Zu}, \cite%
{FZ}, \cite{FI}, \cite{Mi} and \cite{Mi2}).

\bigskip

Observability and controllability for the one-dimensional degenerate
parabolic operator has been extensively studied in many ways: Backstepping
approach for closed-loop control (see \cite{GLM} and \cite{LM}); Carleman
inequalities (see \cite{ABCF}, \cite{CMV}, \cite{CMV2} and \cite{CTY});
Flatness approach (see \cite{Mo} and \cite{BLR}); Moment method (see \cite%
{CMV3} and \cite{CMV4}).

\bigskip

We shall consider the linear unbounded operators ${\mathcal{P}}$ in $%
L^{2}\left( 0,1\right) $, defined by 
\begin{equation*}
\left\{ 
\begin{array}{ll}
{\mathcal{P}}=-\frac{d}{dx}\left( x^{\alpha }\frac{d}{dx}\right) \text{ ,
with }\alpha \in \left[ 0,2\right) \text{ ,} &  \\ 
D({\mathcal{P}})=\left\{ \vartheta \in H_{\alpha ,0}^{1}\left( 0,1\right) 
\text{; }{\mathcal{P}}\vartheta \in L^{2}(0,1)\text{ and BC}_{\alpha
}(\vartheta )=0\right\} \text{ ,} & 
\end{array}%
\right.
\end{equation*}%
where%
\begin{equation*}
H_{\alpha ,0}^{1}(0,1):=\left\{ \vartheta \in L^{2}(0,1)\text{;}\;\vartheta 
\text{ is loc. absolutely continuous in }\left( 0,1\right] \text{,}%
\;\int_{0}^{1}x^{\alpha }|\vartheta ^{\prime }|^{2}<\infty \text{, }%
\vartheta (1)=0\right\} \text{ ,}
\end{equation*}%
and 
\begin{equation*}
\text{BC}_{\alpha }(\vartheta )=%
\begin{cases}
\vartheta _{|_{x=0}}\text{ ,} & \text{for }\alpha \in \lbrack 0,1)\text{ ,}
\\ 
(x^{\alpha }\vartheta ^{\prime })_{|_{x=0}}\text{ ,} & \text{for }\alpha \in
\lbrack 1,2)\text{ .}%
\end{cases}%
\end{equation*}%
Such ${\mathcal{P}}$ is a closed self-adjoint positive densely defined
operator, with compact resolvent. As a consequence, the following spectral
decomposition holds: There exists a countable family of eigenfunctions $\Phi
_{j}$ associated with eigenvalues $\lambda _{j}$ such that

\begin{itemize}
\item $\left\{ \Phi _{j}\right\} _{j\geq 1}$ forms a Hilbert basis of $%
L^{2}(0,1)$

\item ${\mathcal{P}}\Phi_{j}=\lambda_{j}\Phi_{j}$

\item $0<\lambda _{1}\leq \lambda _{2}\leq \dots \leq \lambda
_{k}\rightarrow +\infty $ .
\end{itemize}

An explicit expression of the eigenvalues is given in \cite{Gu} for the
weakly degenerate case $\alpha \in \left( 0,1\right) $, and in \cite{Mo} for
the strongly degenerate case $\alpha \in \lbrack 1,2)$, and depends on the
Bessel functions of first kind (see \cite{MM}). The eigenvalues are simples
and more proprieties are emphasized by Cannarsa, Martinez and Vancostenoble:
First, a uniform bound for the first eigenvalue 
\begin{equation}
\exists c_{1},c_{2}>0\text{ \quad }\forall \alpha \in \lbrack 0,2)\text{
\quad }c_{1}\leq \lambda _{1}\leq c_{2}  \tag{1.2}  \label{1.2}
\end{equation}%
(see \cite{CMV3} (10) at page 176 and (34) at page 183 for $\alpha \in
\lbrack 0,1)$; see \cite{CMV4} proposition 2.13 at page 10 and (3.8)-(3.9)
at page 13 for $\alpha \in \lbrack 1,2)$); Secondly, a uniform spectral gap 
\begin{equation}
\exists \gamma >0\text{ }\quad \forall \alpha \in \lbrack 0,2)\text{ }\quad
\forall k\geq 1\quad \sqrt{\lambda _{k+1}}-\sqrt{\lambda _{k}}\geq \gamma
(2-\alpha )  \tag{1.3}  \label{1.3}
\end{equation}%
(see \cite{CMV3} (74) at page 198 for $\alpha \in \lbrack 0,1)$; see \cite%
{CMV4} at page 30 for $\alpha \in \lbrack 1,2)$).

\bigskip

We are interested in the spectral inequality for the sum of eigenfunctions.
Such Lebeau-Robbiano estimate is done with explicit dependence on $\alpha
\in \lbrack 0,2)$. Our main result is as follows.

\bigskip

Theorem 1.1 .- \textit{Let }$\omega $\textit{\ be an open and nonempty
subset of }$(0,1)$\textit{. There exists a constant }$C>0$\textit{\ such
that }%
\begin{equation*}
\sum_{\lambda _{j}\leq \Lambda }|a_{j}|^{2}\leq Ce^{C\frac{1}{\left(
2-\alpha \right) ^{2}}\sqrt{\Lambda }}\int_{\omega }\left\vert \sum_{\lambda
_{j}\leq \Lambda }a_{j}\Phi _{j}\right\vert ^{2}\text{ ,}
\end{equation*}%
\textit{for any }$\alpha \in \lbrack 0,2)$, $\left\{ a_{j}\right\} \in 
\mathbb{R}$\textit{\ and any }$\Lambda >0$\textit{. }

\bigskip

This is equivalent to 
\begin{equation*}
\sum_{j=1,\cdot \cdot ,N}|a_{j}|^{2}\leq Ce^{C\frac{1}{\left( 2-\alpha
\right) ^{2}}\sqrt{\lambda _{N}}}\int_{\omega }\left\vert \sum_{j=1,\cdot
\cdot ,N}a_{j}\Phi _{j}\right\vert ^{2}\text{ ,}
\end{equation*}%
for any\textit{\ }$\alpha \in \lbrack 0,2)$, $\left\{ a_{j}\right\} \in 
\mathbb{R}$ and any $N>0$.

\bigskip

Here, our approach is based on a combinaison of both Carleman techniques and the
moment method for an elliptic equation. In one hand, it seems difficult to
find the appropriate weight function in Carleman techniques or logarithmic convexity methods for getting
directly the desired spectral inequality. On the other hand, the moment
method is an appropriate tool to get the cost of controllability for the
one-dimensional degenerate parabolic operator.  

\bigskip

As a consequence of Theorem 1.1, we have the following observability
estimate from a measurable set in time for the one-dimensional
degenerate parabolic operator.

\bigskip

Theorem 1.2 .- \textit{Let }$\omega $\textit{\ be an open and nonempty
subset of }$(0,1)$\textit{\ and }$E\subset \left( 0,T\right) $\textit{\ be a
measurable set of positive measure. There exists a constant }$C>0$\textit{\
such that }%
\begin{equation*}
\left\Vert e^{-T{\mathcal{P}}}y_{0}\right\Vert _{L^{2}\left( 0,1\right)
}\leq Ce^{C\frac{1}{(2-\alpha )^{4}}}\int_{\omega \times E}\left\vert e^{-t{%
\mathcal{P}}}y_{0}\right\vert \text{ ,}
\end{equation*}%
\textit{for any }$\alpha \in \lbrack 0,2)$\textit{\ and any }$y_{0}\in
L^{2}\left( 0,1\right) $\textit{. }

\bigskip

This is equivalent to 
\begin{equation*}
\left\Vert y\left( \cdot ,T\right) \right\Vert _{L^{2}\left( 0,1\right)
}\leq Ce^{C\frac{1}{(2-\alpha )^{4}}}\int_{\omega \times E}\left\vert
y\left( x,t\right) \right\vert dxdt\text{ ,}
\end{equation*}%
for any $\alpha \in \lbrack 0,2)$ and any $y_{0}\in L^{2}\left( 0,1\right) $%
\ where $y$\ is the weak solution of the degenerate heat equation%
\begin{equation*}
\left\{ 
\begin{array}{ll}
\partial _{t}y-\partial _{x}\left( x^{\alpha }\partial _{x}y\right) =0\text{
,} & \text{ in }(0,1)\times \left( 0,T\right) \text{ ,} \\ 
\text{BC}_{\alpha }(y)=0\text{ ,} & \text{ on }\{0\}\times \left( 0,T\right) 
\text{ ,} \\ 
y\left( 1,t\right) =0\text{ ,} & \text{ }t\in \left( 0,T\right) \text{ ,} \\ 
y\left( x,0\right) =y_{0}\text{ ,} & \text{ }x\in \left( 0,1\right) \text{ .}%
\end{array}%
\right.
\end{equation*}

\bigskip

In the last years, a lot of works have been devoted to the observability on
measurable sets (see e.g. \cite{AE}, \cite{EMZ}, \cite{PW}, \cite{WZ}, \cite%
{LiZ}). Applications to impulse control and finite-time stabilization can be
rewritten as in \cite{BP}.

\bigskip

\section{Elliptic observation estimates (proof of Theorem 1.1)}

\bigskip

In this Section, our aim is to prove Theorem 1.1 and we start with
presenting the following three results: We first have an uniform
observability estimate for a single eigenfunction given by Proposition 2.1;
Proposition 2.2 establishes a quantitative Holder type of estimate for an
elliptic equation far from the degeneracy; Proposition 2.3 is an uniform
observability estimate for the elliptic equation; We end this Section with
the proof of Theorem 1.1.

\bigskip

Proposition 2.1 .- \textit{For any} $\omega $ \textit{open and nonempty
subset of} $\left( 0,1\right) $, 
\begin{equation*}
\exists \rho >0\text{ \quad }\forall \alpha \in \lbrack 0,2)\text{ \quad }%
\forall j\geq 1\text{ \quad }\int_{\omega }\left\vert \Phi _{j}\right\vert
^{2}\geq \rho (2-\alpha )\text{ .}
\end{equation*}

\bigskip

Given $T>0$ arbitrary, we now consider the following homogeneous elliptic
problem: 
\begin{equation}
\left\{ 
\begin{array}{ll}
\partial _{t}^{2}\varphi -{\mathcal{P}}\varphi =0\text{ ,} & \text{in }%
\left( 0,1\right) \times \left( 0,T\right) \text{ ,} \\ 
\text{BC}_{\alpha }(\varphi )=0\text{ ,} & \text{on }\{0\}\times \left(
0,T\right) \text{ ,} \\ 
\varphi _{|_{x=1}}=0\text{ ,} & \text{on }\{1\}\times \left( 0,T\right) 
\text{ ,} \\ 
\varphi \left( \cdot ,0\right) =\varphi _{0}\text{ ,} & \text{in }\left(
0,1\right) \text{ ,} \\ 
\partial _{t}\varphi \left( \cdot ,0\right) =\varphi _{1}\text{ ,} & \text{%
in }\left( 0,1\right) \text{ ,}%
\end{array}%
\right.  \tag{2.1}  \label{2.1}
\end{equation}%
where $\varphi _{0}$ and $\varphi _{1}$ belong to span$\left\{ \Phi
_{j};1\leq j\leq N\right\} $.

\bigskip

Proposition 2.2 .- \textit{Let }$0<a<b<1$\textit{\ and }$T>0$\textit{. There
exist }$c>0$\textit{\ and }$\delta \in \left( 0,1\right) $\textit{\ such
that for all }$\alpha \in \left[ 0,2\right) $\textit{, the solution }$%
\varphi $ \textit{of (\ref{2.1})\ satisfies}%
\begin{equation*}
\left\Vert \varphi \right\Vert _{H^{1}\left( \left( \frac{2a+b}{3},\frac{a+2b%
}{3}\right) \times \left( 0,T/4\right) \right) }\leq c\left\Vert \varphi
\right\Vert _{H^{1}\left( \left( a.b\right) \times \left( 0,T\right) \right)
}^{1-\delta }\left( \left\Vert \varphi _{0}\right\Vert _{H^{1}\left(
a,b\right) }+\left\Vert \varphi _{1}\right\Vert _{L^{2}\left( a,b\right)
}\right) ^{\delta }\text{ ,}
\end{equation*}

\bigskip

Proposition 2.3 .- \textit{Let} $\omega $ \textit{be an open and nonempty
subset of} $\left( 0,1\right) $. \textit{For any} $N\geq 1$\textit{, }$T>0$%
\textit{, and any }$\alpha \in \lbrack 0,2)$\textit{, the solution }$\varphi 
$\textit{\ of (\ref{2.1}) satisfies} 
\begin{equation*}
\left\Vert \varphi \left( \cdot ,T\right) \right\Vert _{L^{2}\left(
0,1\right) }^{2}\leq \frac{C\left( 1+\lambda _{N}\right) }{\rho
^{2}(2-\alpha )^{2}}\left( 1+\frac{1}{T}\right) e^{C\sqrt{\lambda _{N}}%
\left( T+\frac{1}{T\gamma ^{2}(2-\alpha )^{2}}\right)
}\int_{0}^{T}\int_{\omega }\left\vert \varphi \right\vert ^{2}\text{ ,}
\end{equation*}%
\textit{where }$C>0$\textit{\ is independent of }$N,T>0$\textit{\ and }$%
\alpha \in \lbrack 0,2)$\textit{. Here }$\rho $\textit{\ is given by
Proposition 2.1 and }$\gamma $\textit{\ comes from (\ref{1.3}).}

\bigskip

Now, we are able to present the proof of Theorem 1.1.

\bigskip

Our strategy is as follows. We will use Proposition 2.3 and we observe the
whole domain (including the region where the ellipticity degenerates) from
one region where the operator $\partial _{t}^{2}+\mathcal{P}$ is uniformly
elliptic; there, we use classical global Carleman techniques to observe from
the boundary $(a,b)\times \{0\}$ with Proposition 2.2. That observation
region provides precisely the right hand side of Theorem 1.1.

\bigskip

Proof of Theorem 1.1 .- We consider the above homogeneous elliptic problem
with $\varphi _{0}\left( x\right) =0$ and $\varphi _{1}\left( x\right)
=\sum\limits_{j=1,\cdot \cdot ,N}a_{j}\Phi _{j}\left( x\right) $ where $%
\left\{ a_{j}\right\} \in \mathbb{R}$. Recall that $\varphi $ can be
explicitly written by Fourier series: For any $x\in \left( 0,1\right) $, 
\begin{equation*}
\varphi \left( x,t\right) =\sum\limits_{j=1,\cdot \cdot ,N}\frac{1}{\sqrt{%
\lambda _{j}}}\text{sinh}\left( \sqrt{\lambda _{j}}t\right) a_{j}\Phi
_{j}\left( x\right) \text{ .}
\end{equation*}%
Let $0<a<b<1$ and set $\omega =\left( a,b\right) $ and $\widetilde{\omega }%
=\left( \frac{2a+b}{3},\frac{a+2b}{3}\right) $. We have , for some constants 
$C,C_{1},C_{2}>0$ independent on $N$ and $\alpha $, 
\begin{equation*}
\begin{array}{lll}
\displaystyle\sum\limits_{j=1,\cdot \cdot ,N}|a_{j}|^{2} & \leq \displaystyle%
\sum\limits_{j=1,\cdot \cdot ,N}|a_{j}|^{2}\frac{1}{\lambda _{j}}\text{sinh}%
^{2}\left( \sqrt{\lambda _{j}}/4\right) Ce^{C\sqrt{\lambda _{N}}} & \text{by
(\ref{1.2})} \\ 
& =Ce^{C\sqrt{\lambda _{N}}}\left\Vert \varphi \left( \cdot ,1/4\right)
\right\Vert _{L^{2}\left( 0,1\right) }^{2} &  \\ 
& \leq C_{1}e^{C_{1}\frac{1}{(2-\alpha )^{2}}\sqrt{\lambda _{N}}}%
\displaystyle\int_{0}^{1/4}\int_{\widetilde{\omega }}\left\vert \varphi
\right\vert ^{2} & \text{by Proposition 2.3 applied to }\widetilde{\omega }%
\times \left( 0,1/4\right) \\ 
& \leq C_{2}e^{C_{2}\frac{1}{(2-\alpha )^{2}}\sqrt{\lambda _{N}}}\left\Vert
\varphi \right\Vert _{H^{1}\left( \omega \times \left( 0,1\right) \right)
}^{2\left( 1-\delta \right) }\left\Vert \varphi _{1}\right\Vert
_{L^{2}\left( \omega \right) }^{2\delta } & \text{by Proposition 2.2.}%
\end{array}%
\end{equation*}%
But, 
\begin{equation*}
\left\Vert \varphi _{1}\right\Vert _{L^{2}\left( \omega \right)
}^{2}=\int_{\omega }\left\vert \sum\limits_{j=1,\cdot \cdot ,N}a_{j}\Phi
_{j}\right\vert ^{2}
\end{equation*}%
and for some constants $c_{1},c_{2}>0$ independent on $N$ and $\alpha $, it
holds 
\begin{eqnarray*}
\left\Vert \varphi \right\Vert _{H^{1}\left( \omega \times \left( 0,1\right)
\right) }^{2} &=&\left\Vert \varphi \right\Vert _{L^{2}\left( \omega \times
\left( 0,1\right) \right) }^{2}+\left\Vert \partial _{x}\varphi \right\Vert
_{L^{2}\left( \omega \times \left( 0,1\right) \right) }^{2}+\left\Vert
\partial _{t}\varphi \right\Vert _{L^{2}\left( \omega \times \left(
0,1\right) \right) }^{2} \\
&\leq &\left\Vert \varphi \right\Vert _{L^{2}\left( \left( 0,1\right)
^{2}\right) }^{2}+c_{1}\left\Vert x^{\alpha /2}\partial _{x}\varphi
\right\Vert _{L^{2}\left( \left( 0,1\right) ^{2}\right) }^{2}+\left\Vert
\partial _{t}\varphi \right\Vert _{L^{2}\left( \left( 0,1\right) ^{2}\right)
}^{2} \\
&\leq &c_{2}e^{c_{2}\sqrt{\lambda _{N}}}\sum\limits_{j=1,\cdot \cdot
,N}|a_{j}|^{2}\text{ .}
\end{eqnarray*}%
Combining the above estimates completes the proof of Theorem 1.1.

\bigskip

\section{Elliptic observability by the moment method (proof of Proposition
2.3)}

\bigskip

In this section, we shall prove Proposition 2.3. Let $\omega $ be an open
and nonempty subset of $\left( 0,1\right) $. Given $T>0$ arbitrary, we
consider the following non-homogeneous elliptic problem: 
\begin{equation}
\left\{ 
\begin{array}{ll}
\partial _{t}^{2}u-{\mathcal{P}}u=h\text{ ,} & \text{in }\left( 0,1\right)
\times \left( 0,T\right) \text{ ,} \\ 
\text{BC}_{\alpha }(u)=0\text{ ,} & \text{on }\{0\}\times \left( 0,T\right) 
\text{ ,} \\ 
u_{|_{x=1}}=0\text{ ,} & \text{on }\{1\}\times \left( 0,T\right) \text{ ,}
\\ 
u\left( \cdot ,0\right) =u_{0}\text{ ,} & \text{in }\left( 0,1\right) \text{
,} \\ 
\partial _{t}u\left( \cdot ,0\right) =u_{1}\text{ ,} & \text{in }\left(
0,1\right) \text{ ,}%
\end{array}%
\right.  \tag{3.1}  \label{3.1}
\end{equation}%
where%
\begin{equation}
\left\{ 
\begin{array}{l}
h\left( x,t\right) =\sum\limits_{j=1,\cdot \cdot ,N}\sum\limits_{k=1,\cdot
\cdot ,N}g_{k}\left( t\right) \left( \int_{\omega }\Phi _{j}\Phi _{k}\right)
\Phi _{j}\left( x\right) \text{ with }g\left( x,t\right)
=\sum\limits_{k=1,\cdot \cdot ,N}g_{k}\left( t\right) \Phi _{k}\left(
x\right) \text{ ,} \\ 
u_{0}\left( x\right) =\sum\limits_{j=1,\cdot \cdot ,N}a_{j}\Phi _{j}\left(
x\right) \text{ ,} \\ 
u_{1}\left( x\right) =\sum\limits_{j=1,\cdot \cdot ,N}b_{j}\Phi _{j}\left(
x\right) \text{ .}%
\end{array}%
\right.  \tag{3.2}  \label{3.2}
\end{equation}

\bigskip

\subsection{Well-posedness property}

\bigskip

Definition 3.1 .- Let $N\in \mathbb{N}^{\ast }$. We denote $\Pi _{N}L^{2}=$%
span$\left\{ \Phi _{j};1\leq j\leq N\right\} $. The space $\Pi _{N}L^{2}$
endowed with the $L^{2}\left( \Omega \right) $ norm is a finite dimensional
Hilbert space.\newline

%Consider the operator
%\[
%A=\begin{pmatrix} 0 & -I \\ - \mathcal P & 0 \end{pmatrix},
%\]
%acting on $(\Pi_{\Lambda}L^{2})^{2}$. As $A$ is bounded, it is the infinitesimal generator of a uniformly continuous semigroup (see Pazy, Theorem 1.2).\\

It is well-known that when $g_{j}\in L^{2}(0,T)$, the unique solution of (%
\ref{3.1}) verifies $u\in H^{2}(0,T;\Pi _{N}L^{2})$ and is given by the
Duhamel formula 
\begin{multline*}
u(\cdot ,t)=\sum\limits_{j=1,\cdot \cdot ,N}\text{cosh}(\sqrt{\lambda _{j}}%
t)a_{j}\Phi _{j}+\sum\limits_{j=1,\cdot \cdot ,N}\frac{\text{sinh}(\sqrt{%
\lambda _{j}}t)}{\sqrt{\lambda _{j}}}b_{j}\Phi _{j} \\
+\sum\limits_{j=1,\cdot \cdot ,N}\sum\limits_{k=1,\cdot \cdot
,N}(\int_{\omega }\Phi _{j}\Phi _{k})\int_{0}^{t}\frac{\text{sinh}(\sqrt{%
\lambda _{j}}(t-s))}{\sqrt{\lambda _{j}}}g_{k}(s)ds\Phi _{j}\text{ .}
\end{multline*}%
%
%
%
%
%
%
%
%
%
%
%
%
%
%
%
%
%
%
%
%
%
%
%
%
%
%
%
%
%
%
%
%
%
%
%The proof follows from [Pazy, Section 4, Theorem 2.9]. Moreover, $\partial_{t}^{2}u = \mathcal P u +h \in L^{2}(0,T;\Pi_{N}L^{2})$.

\bigskip

\subsection{Construction of the control}

\bigskip

Definition 3.2 .- We say that system (\ref{3.1}) is controllable at time $T$
if for any $\left( u_{0},u_{1}\right) \in (\Pi _{N}L^{2})^{2}$ there is $%
g\in L^{2}(0,T;\Pi _{N}L^{2})$ as in (\ref{3.2}) such that 
\begin{equation*}
u\left( \cdot ,T\right) =\partial _{t}u\left( \cdot ,T\right) =0\text{ .}
\end{equation*}

\bigskip

Lemma 3.1 .- \textit{Equation (\ref{3.1}) is controllable in time }$T$%
\textit{\ if and only if, for any }$\left( u_{0},u_{1}\right) \in (\Pi
_{N}L^{2})^{2}$ \textit{there is} $g\in L^{2}(0,T;\Pi _{N}L^{2})$\textit{\
as in (\ref{3.2}) such that the following relation holds} 
\begin{equation}
-\int_{0}^{1}u_{1}\varphi \left( \cdot ,T\right) -\int_{0}^{1}u_{0}\partial
_{t}\varphi \left( \cdot ,T\right) =\int_{0}^{T}\int_{\omega }g\left(
x,t\right) \varphi \left( x,T-t\right) dxdt  \tag{3.3}  \label{3.3}
\end{equation}%
\textit{for any }$\left( \varphi _{0},\varphi _{1}\right) \in (\Pi
_{N}L^{2})^{2}$\textit{, where }$\varphi $\textit{\ is the solution of (\ref%
{2.1}).}

\textit{Further, if the system (\ref{3.1}) is controllable at time }$T$%
\textit{\ with a control }$g\in L^{2}(0,T;\Pi _{N}L^{2})$\textit{\
satisfying the bound} 
\begin{equation*}
\left\Vert g\right\Vert _{L^{2}\left( \left( 0,1\right) \times \left(
0,T\right) \right) }^{2}:=\sum\limits_{j=1,\cdot \cdot
,N}\int_{0}^{T}\left\vert g_{j}\left( t\right) \right\vert ^{2}\leq
K\left\Vert \left( u_{0},u_{1}\right) \right\Vert _{\left( L^{2}\left(
0,1\right) \right) ^{2}}^{2}:=K\sum\limits_{j=1,\cdot \cdot ,N}\left(
a_{j}^{2}+b_{j}^{2}\right)
\end{equation*}%
\textit{for some }$K>0$\textit{, then the solution }$\varphi $\textit{\ of (%
\ref{2.1}) satisfies} 
\begin{equation*}
\left\Vert \varphi \left( \cdot ,T\right) \right\Vert _{L^{2}\left(
0,1\right) }^{2}+\left\Vert \partial _{t}\varphi \left( \cdot ,T\right)
\right\Vert _{L^{2}\left( 0,1\right) }^{2}\leq K\int_{0}^{T}\int_{\omega
}\left\vert \varphi \right\vert ^{2}\text{ .}
\end{equation*}

\bigskip

Proof of Lemma 3.1 .- Let $g\in L^{2}(0,T;\Pi _{N}L^{2})$ be arbitrary and $%
u $ be the solution of (\ref{3.1}). Given $\varphi $ the solution of (\ref%
{2.1}) then, by multiplying (\ref{3.1}) by $\varphi \left( x,T-t\right) $
and by integrating by parts we obtain that%
\begin{equation*}
\int_{0}^{1}\partial _{t}u\left( \cdot ,T\right) \varphi
_{0}+\int_{0}^{1}u\left( \cdot ,T\right) \varphi
_{1}-\int_{0}^{1}u_{1}\varphi \left( \cdot ,T\right)
-\int_{0}^{1}u_{0}\partial _{t}\varphi \left( \cdot ,T\right)
=\int_{0}^{T}\int_{0}^{1}h\left( x,t\right) \varphi \left( x,T-t\right) dxdt%
\text{ }
\end{equation*}%
and 
\begin{equation*}
\int_{0}^{T}\int_{0}^{1}h\left( x,t\right) \varphi \left( x,T-t\right)
dxdt=\int_{0}^{T}\int_{\omega }g\left( x,t\right) \varphi \left(
x,T-t\right) dxdt\text{ .}
\end{equation*}%
Now, if (\ref{3.3}) is verified, it follows that 
\begin{equation*}
\int_{0}^{1}\partial _{t}u\left( \cdot ,T\right) \varphi
_{0}+\int_{0}^{1}u\left( \cdot ,T\right) \varphi _{1}=0
\end{equation*}%
for any $\left( \varphi _{0},\varphi _{1}\right) \in (\Pi _{N}L^{2})^{2}$
which implies that $u\left( \cdot ,T\right) =\partial _{t}u\left( \cdot
,T\right) =0$. Hence, the solution is controllable at time $T$ and $g$ is a
control for (\ref{3.1}). Reciprocally, if $g\in L^{2}(0,T;\Pi _{N}L^{2})$ is
a control for (\ref{3.1}), we have that $u\left( \cdot ,T\right) =\partial
_{t}u\left( \cdot ,T\right) =0$. It implies that (\ref{3.3}) holds. Finally,
one can choose $\left( u_{0},u_{1}\right) =\left( \partial _{t}\varphi
\left( \cdot ,T\right) ,\varphi \left( \cdot ,T\right) \right) $ and apply (%
\ref{3.3}) to get the desired estimate thanks to Cauchy-Schwarz inequality
and the proof finishes.

\bigskip

Proof of Proposition 2.3 .- Our aim is to construct a control $g$ given by $%
g\left( x,t\right) =\sum\limits_{k=1,\cdot \cdot ,N}g_{k}\left( t\right)
\Phi _{k}\left( x\right) $ such that (\ref{3.3}) holds. Let 
\begin{equation*}
\left\{ 
\begin{array}{l}
\varphi _{0}\left( x\right) =\sum\limits_{j=1,\cdot \cdot ,N}c_{j}\Phi
_{j}\left( x\right) \text{ ,} \\ 
\varphi _{1}\left( x\right) =\sum\limits_{j=1,\cdot \cdot ,N}d_{j}\Phi
_{j}\left( x\right) \text{ }%
\end{array}%
\right.
\end{equation*}%
be the initial data of (\ref{2.1}). Then, recall that $\varphi $ can be
explicitly written by Fourier series: For any $x\in \left( 0,1\right) $, 
\begin{equation*}
\varphi \left( x,t\right) =\sum\limits_{j=1,\cdot \cdot ,N}\left( e^{\sqrt{%
\lambda _{j}}t}\frac{1}{2}\left( c_{j}+\frac{1}{\sqrt{\lambda _{j}}}%
d_{j}\right) +e^{-\sqrt{\lambda _{j}}t}\frac{1}{2}\left( c_{j}-\frac{1}{%
\sqrt{\lambda _{j}}}d_{j}\right) \right) \Phi _{j}\left( x\right) \text{ .}
\end{equation*}

\bigskip

First, let us clarify the expression $-\int_{0}^{1}u_{1}\varphi \left( \cdot
,T\right) -\int_{0}^{1}u_{0}\partial _{t}\varphi \left( \cdot ,T\right) $ :%
\begin{equation}
\begin{array}{ll}
-\displaystyle\int_{0}^{1}u_{1}\varphi \left( \cdot ,T\right) -\displaystyle%
\int_{0}^{1}u_{0}\partial _{t}\varphi \left( \cdot ,T\right)  & =%
\displaystyle\sum\limits_{j=1,\cdot \cdot ,N}e^{\sqrt{\lambda _{j}}T}\frac{1%
}{2}\left( c_{j}+\frac{1}{\sqrt{\lambda _{j}}}d_{j}\right) \displaystyle%
\int_{0}^{1}\left( -u_{1}-\sqrt{\lambda _{j}}u_{0}\right) \Phi _{j} \\ 
& \text{\quad }+\displaystyle\sum\limits_{j=1,\cdot \cdot ,N}e^{-\sqrt{%
\lambda _{j}}T}\frac{1}{2}\left( c_{j}-\frac{1}{\sqrt{\lambda _{j}}}%
d_{j}\right) \displaystyle\int_{0}^{1}\left( -u_{1}+\sqrt{\lambda _{j}}%
u_{0}\right) \Phi _{j}\text{ .}%
\end{array}
\tag{3.4}  \label{3.4}
\end{equation}

Next, let us clarify the expression $\int_{0}^{T}\int_{\omega }g\left(
x,t\right) \varphi \left( x,T-t\right) dxdt$, that is $\int_{0}^{T}%
\int_{0}^{1}h\left( x,t\right) \varphi \left( x,T-t\right) dxdt$ :%
\begin{eqnarray*}
&&\int_{0}^{T}\int_{\omega }g\left( x,t\right) \varphi \left( x,T-t\right)
dxdt \\
&=&\sum\limits_{k=1,\cdot \cdot ,N}\sum\limits_{j=1,\cdot \cdot ,N}e^{\sqrt{%
\lambda _{j}}T}\frac{1}{2}\left( c_{j}+\frac{1}{\sqrt{\lambda _{j}}}%
d_{j}\right) \int_{\omega }\Phi _{k}\Phi _{j}\int_{0}^{T}g_{k}\left(
t\right) e^{-\sqrt{\lambda _{j}}t}dt \\
&&+\sum\limits_{k=1,\cdot \cdot ,N}\sum\limits_{j=1,\cdot \cdot ,N}e^{-\sqrt{%
\lambda _{j}}T}\frac{1}{2}\left( c_{j}-\frac{1}{\sqrt{\lambda _{j}}}%
d_{j}\right) \int_{\omega }\Phi _{k}\Phi _{j}\int_{0}^{T}g_{k}\left(
t\right) e^{\sqrt{\lambda _{j}}t}dt\text{ .}
\end{eqnarray*}%
Now, suppose that $g_{k}\left( t\right) =\alpha _{k}\sigma _{k}^{0}\left(
t\right) +\beta _{k}\sigma _{k}^{1}\left( t\right) $ where $\sigma _{k}^{0}$%
, $\sigma _{k}^{1}$ belong to $L^{2}\left( 0,T\right) $ and that the
following moment formula holds: 
\begin{equation}
\left\{ 
\begin{array}{c}
\displaystyle\int_{0}^{T}\sigma _{k}^{0}\left( t\right) e^{-\sqrt{\lambda
_{j}}t}dt=0\text{ and }\displaystyle\int_{0}^{T}\sigma _{k}^{1}\left(
t\right) e^{-\sqrt{\lambda _{j}}t}dt=\delta _{jk}\text{ ;} \\ 
\displaystyle\int_{0}^{T}\sigma _{k}^{0}\left( t\right) e^{\sqrt{\lambda _{j}%
}t}dt=\delta _{jk}\text{ and }\int_{0}^{\displaystyle T}\sigma
_{k}^{1}\left( t\right) e^{\sqrt{\lambda _{j}}t}dt=0\text{ ,}%
\end{array}%
\right.   \tag{3.5}  \label{3.5}
\end{equation}%
then, we obtain%
\begin{equation}
\begin{array}{ll}
\displaystyle\int_{0}^{T}\int_{\omega }g\left( x,t\right) \varphi \left(
x,T-t\right) dxdt & =\displaystyle\sum\limits_{j=1,\cdot \cdot ,N}e^{\sqrt{%
\lambda _{j}}T}\frac{1}{2}\left( c_{j}+\frac{1}{\sqrt{\lambda _{j}}}%
d_{j}\right) \beta _{j}\displaystyle\int_{\omega }\left\vert \Phi
_{j}\right\vert ^{2} \\ 
& +\displaystyle\sum\limits_{j=1,\cdot \cdot ,N}e^{-\sqrt{\lambda _{j}}T}%
\frac{1}{2}\left( c_{j}-\frac{1}{\sqrt{\lambda _{j}}}d_{j}\right) \alpha _{j}%
\displaystyle\int_{\omega }\left\vert \Phi _{j}\right\vert ^{2}\text{ .}%
\end{array}
\tag{3.6}  \label{3.6}
\end{equation}
By comparing the identities (\ref{3.4}) and (\ref{3.6}), one can deduce that
if 
\begin{equation*}
\int_{0}^{1}\left( -u_{1}-\sqrt{\lambda _{j}}u_{0}\right) \Phi _{j}=\beta
_{j}\int_{\omega }\left\vert \Phi _{j}\right\vert ^{2}\text{ and }%
\int_{0}^{1}\left( -u_{1}+\sqrt{\lambda _{j}}u_{0}\right) \Phi _{j}=\alpha
_{j}\int_{\omega }\left\vert \Phi _{j}\right\vert ^{2}
\end{equation*}%
for any $j=1,\cdot \cdot ,N$, then (\ref{3.3}) holds for any $\left( \varphi
_{0},\varphi _{1}\right) $ which implies by Lemma 3.1 that (\ref{3.1}) is
controllable in time $T$.

\bigskip

Therefore, one can conclude that the control given by $g\left( x,t\right)
:=\sum\limits_{j=1,\cdot \cdot ,N}\left[ \alpha _{j}\sigma _{j}^{0}\left(
t\right) +\beta _{j}\sigma _{j}^{1}\left( t\right) \right] \Phi _{j}\left(
x\right) $ where%
\begin{equation*}
\alpha _{j}:=\frac{\int_{0}^{1}\left( -u_{1}+\sqrt{\lambda _{j}}u_{0}\right)
\Phi _{j}}{\int_{\omega }\left\vert \Phi _{j}\right\vert ^{2}}=\frac{-b_{j}+%
\sqrt{\lambda _{j}}a_{j}}{\int_{\omega }\left\vert \Phi _{j}\right\vert ^{2}}%
\text{ and }\beta _{j}:=\frac{\int_{0}^{1}\left( -u_{1}-\sqrt{\lambda _{j}}%
u_{0}\right) \Phi _{j}}{\int_{\omega }\left\vert \Phi _{j}\right\vert ^{2}}=%
\frac{-b_{j}-\sqrt{\lambda _{j}}a_{j}}{\int_{\omega }\left\vert \Phi
_{j}\right\vert ^{2}}\text{ ,}
\end{equation*}%
is an appropriate candidate. Notice that by Proposition 2.1, $\int_{\omega
}\left\vert \Phi _{j}\right\vert ^{2}\neq 0$. It remains to construct the
sequence of functions $\left( \sigma _{k}^{0},\sigma _{k}^{1}\right) _{k\geq
1}$ in $\left( L^{2}\left( 0,T\right) \right) ^{2}$ such that (\ref{3.5})
holds. Such property is called biorthogonality of the family $\left( \sigma
_{k}^{0},\sigma _{k}^{1}\right) _{k\geq 1}$. To do so, we apply the
following result from Cannarsa, Martinez and Vancostenoble (see \cite{CMV3}
Theorem 2.4 at page 179) :

\bigskip

Theorem 3.1 .- (Existence of a suitable biorthogonal family and upper
bounds) Assume that 
\begin{equation*}
\forall n>0\text{, }\mu _{n}\geq 0
\end{equation*}%
and that there is some $r>0$ such that%
\begin{equation*}
\forall n>0\text{, }\sqrt{\mu _{n+1}}-\sqrt{\mu _{n}}\geq r\text{ .}
\end{equation*}%
Then there exists a family $\left( \theta _{m}\right) _{m>0}$ which is
biorthogonal to the family $\left( e^{\mu _{n}t}\right) _{n>0}$ in $%
L^{2}(0,T)$:%
\begin{equation*}
\forall m,n>0\text{, }\int_{0}^{T}\theta _{m}\left( t\right) e^{\mu
_{n}t}dt=\delta _{mn}\text{ .}
\end{equation*}%
Moreover, it satisfies: there is some universal constant $c$ independent of $%
T$, $r$ and $m$ such that, for all $m>0$, we have%
\begin{equation*}
\left\Vert \theta _{m}\right\Vert _{L^{2}\left( 0,T\right) }^{2}\leq
ce^{-2\mu _{m}T}e^{c\frac{1}{r}\sqrt{\mu _{m}}}B\left( T,r\right)
\end{equation*}%
with%
\begin{equation*}
B\left( T,r\right) =\left\{ 
\begin{array}{cc}
\left( \frac{1}{T}+\frac{1}{T^{2}r^{2}}\right) e^{c\frac{1}{Tr^{2}}} & \text{%
if }T\leq \frac{1}{r^{2}}, \\ 
cr^{2} & \text{if }T\geq \frac{1}{r^{2}}\text{ .}%
\end{array}%
\right.
\end{equation*}

\bigskip

Now, define the increasing sequence of non negative real numbers $\left( \mu
_{n}\right) _{n\geq 1}$ as follows:%
\begin{equation*}
\mu _{n}=\left\vert 
\begin{array}{ll}
\sqrt{\lambda _{N}}-\sqrt{\lambda _{N-\left( n-1\right) }} & \text{if }1\leq
n\leq N ,\\ 
\sqrt{\lambda _{N}}+\sqrt{\lambda _{n-N}} & \text{if }N+1\leq n\leq 2N, \\ 
\left( \sqrt{\mu _{n-1}}+\gamma \left( \lambda _{N}\right) ^{-1/4}\right)
^{2} & \text{if }n\geq 2N+1\text{ .}%
\end{array}%
\right.
\end{equation*}%
We need to check that such sequence fulfills the assumption of Theorem 3.1
thanks to the fact that $\sqrt{\lambda _{k+1}}-\sqrt{\lambda _{k}}\geq
\gamma (2-\alpha )$ given by (\ref{1.3}). Indeed, for any $1\leq n\leq N-1$, 
\begin{equation*}
\sqrt{\mu _{n+1}}-\sqrt{\mu _{n}}=\frac{\sqrt{\lambda _{N-\left( n-1\right) }%
}-\sqrt{\lambda _{N-n}}}{\sqrt{\sqrt{\lambda _{N}}-\sqrt{\lambda _{N-n}}}+%
\sqrt{\sqrt{\lambda _{N}}-\sqrt{\lambda _{N-\left( n-1\right) }}}}\geq \frac{%
\gamma (2-\alpha )}{2\left( \lambda _{N}\right) ^{1/4}}\text{ ;}
\end{equation*}%
for any $N+1\leq n\leq 2N-1$, 
\begin{equation*}
\sqrt{\mu _{n+1}}-\sqrt{\mu _{n}}=\frac{\sqrt{\lambda _{n+1-N}}-\sqrt{%
\lambda _{n-N}}}{\sqrt{\sqrt{\lambda _{N}}+\sqrt{\lambda _{n+1-N}}}+\sqrt{%
\sqrt{\lambda _{N}}+\sqrt{\lambda _{n-N}}}}\geq \frac{\gamma (2-\alpha )}{2%
\sqrt{2}\left( \lambda _{N}\right) ^{1/4}}\text{ ;}
\end{equation*}%
for any $n\geq 2N$, $\sqrt{\mu _{n+1}}-\sqrt{\mu _{n}}=\gamma \left( \lambda
_{N}\right) ^{-1/4}$ and 
\begin{equation*}
\sqrt{\mu _{N+1}}-\sqrt{\mu _{N}}=\frac{2\sqrt{\lambda _{1}}}{\sqrt{\sqrt{%
\lambda _{N}}+\sqrt{\lambda _{1}}}+\sqrt{\sqrt{\lambda _{N}}-\sqrt{\lambda
_{1}}}}\geq \frac{2\sqrt{\lambda _{1}}}{\left( 1+\sqrt{2}\right) \left(
\lambda _{N}\right) ^{1/4}}\text{ .}
\end{equation*}%
Consequently, it fulfills by a straightforward computation the assumptions
of the above Theorem 3.1: Precisely, 
\begin{equation*}
\forall n>0\text{, }\mu _{n}\geq 0\text{ and }\sqrt{\mu _{n+1}}-\sqrt{\mu
_{n}}\geq r\text{ ,}
\end{equation*}%
with 
\begin{equation}
r=\frac{\varsigma }{\left( \lambda _{N}\right) ^{1/4}}\text{ , and\quad }%
\varsigma =\text{min}\left( \frac{\gamma (2-\alpha )}{2\sqrt{2}},\frac{2%
\sqrt{\lambda _{1}}}{1+\sqrt{2}}\right) \text{ .}  \tag{3.7}  \label{3.7}
\end{equation}%
By Theorem 3.1, we have a family $\left( \theta _{m}\right) _{m>0}$ which is
biorthogonal to the family $\left( e^{\mu _{n}t}\right) _{n>0}$ in $%
L^{2}(0,T)$:%
\begin{equation*}
\forall m,n>0\text{, }\int_{0}^{T}\theta _{m}\left( t\right) e^{\mu
_{n}t}dt=\delta _{mn}\text{ .}
\end{equation*}%
Therefore, 
\begin{equation*}
\text{if }1\leq n\leq N\text{, then }\int_{0}^{T}\theta _{m}\left( t\right)
e^{\sqrt{\lambda _{N}}t}e^{-\sqrt{\lambda _{N-\left( n-1\right) }}%
t}dt=\delta _{mn}\text{ ;}
\end{equation*}%
\begin{equation*}
\text{if }N+1\leq n\leq 2N\text{, then }\int_{0}^{T}\theta _{m}\left(
t\right) e^{\sqrt{\lambda _{N}}t}e^{\sqrt{\lambda _{n-N}}t}dt=\delta _{mn}%
\text{ .}
\end{equation*}%
That is, for any $j=1,\cdot \cdot ,N$, 
\begin{equation}
\int_{0}^{T}\theta _{N-\left( j-1\right) }\left( t\right) e^{\sqrt{\lambda
_{N}}t}e^{-\sqrt{\lambda _{j}}t}dt=1\text{ ; }\int_{0}^{T}\theta _{m}\left(
t\right) e^{\sqrt{\lambda _{N}}t}e^{-\sqrt{\lambda _{j}}t}dt=0\text{ when }%
m\neq N-\left( j-1\right) \text{ ;}  \tag{3.8}  \label{3.8}
\end{equation}%
\begin{equation}
\int_{0}^{T}\theta _{N+j}\left( t\right) e^{\sqrt{\lambda _{N}}t}e^{\sqrt{%
\lambda _{j}}t}dt=1\text{ ; }\int_{0}^{T}\theta _{m}\left( t\right) e^{\sqrt{%
\lambda _{N}}t}e^{\sqrt{\lambda _{j}}t}dt=0\text{ when }m\neq N+j\text{ .} 
\tag{3.9}  \label{3.9}
\end{equation}%
Finally, we set for any $k=1,\cdot \cdot ,N$, 
\begin{equation*}
\sigma _{k}^{0}\left( t\right) =\theta _{N+k}\left( t\right) e^{\sqrt{%
\lambda _{N}}t}\text{ and }\sigma _{k}^{1}\left( t\right) =\theta _{N-\left(
k-1\right) }\left( t\right) e^{\sqrt{\lambda _{N}}t}
\end{equation*}%
in order that by (\ref{3.8}), for $k,j=1,\dots ,N,$ 
\begin{equation*}
\int_{0}^{T}\sigma _{k}^{0}\left( t\right) e^{-\sqrt{\lambda _{j}}t}dt=0%
\text{ and }\int_{0}^{T}\sigma _{k}^{1}\left( t\right) e^{-\sqrt{\lambda _{j}%
}t}dt=\delta _{jk}
\end{equation*}%
and by (\ref{3.9})%
\begin{equation*}
\int_{0}^{T}\sigma _{k}^{0}\left( t\right) e^{\sqrt{\lambda _{j}}t}dt=\delta
_{jk}\text{ and }\int_{0}^{T}\sigma _{k}^{1}\left( t\right) e^{\sqrt{\lambda
_{j}}t}dt=0\text{ .}
\end{equation*}%
Further, it holds that for any $k=1,\cdot \cdot ,N$,%
\begin{equation}
\left\Vert \sigma _{k}^{0}\right\Vert _{L^{2}\left( 0,T\right) }^{2}\leq e^{2%
\sqrt{\lambda _{N}}T}\left\Vert \theta _{N+k}\right\Vert _{L^{2}\left(
0,T\right) }^{2}\text{ and }\left\Vert \sigma _{k}^{1}\right\Vert
_{L^{2}\left( 0,T\right) }^{2}\leq e^{2\sqrt{\lambda _{N}}T}\left\Vert
\theta _{N-\left( k-1\right) }\right\Vert _{L^{2}\left( 0,T\right) }^{2}%
\text{ .}  \tag{3.10}  \label{3.10}
\end{equation}%
This completes the construction of our control given by $g\left( x,t\right)
:=\sum\limits_{j=1,\cdot \cdot ,N}\left[ \alpha _{j}\sigma _{j}^{0}\left(
t\right) +\beta _{j}\sigma _{j}^{1}\left( t\right) \right] \Phi _{j}\left(
x\right) $.

\bigskip

\subsection{Cost of the control}

\bigskip

Theorem 3.1 with (\ref{3.7}) implies that there is some universal constant $%
c $ independent of $T$ and $N$ such that for any $m=1,\cdot \cdot ,2N$, 
\begin{eqnarray*}
\left\Vert \theta _{m}\right\Vert _{L^{2}\left( 0,T\right) }^{2} &\leq &ce^{c%
\frac{1}{r}\sqrt{\mu _{m}}}B\left( T,r\right) :=ce^{c\frac{\left( \lambda
_{N}\right) ^{1/4}}{\varsigma }\sqrt{\mu _{m}}}B\left( T,\varsigma \left(
\lambda _{N}\right) ^{-1/4}\right) \\
&\leq &ce^{\frac{c\sqrt{2}}{\varsigma }\sqrt{\lambda _{N}}}B\left(
T,\varsigma \left( \lambda _{N}\right) ^{-1/4}\right)
\end{eqnarray*}%
because $\sqrt{\mu _{m}}\leq \sqrt{2}\left( \lambda _{N}\right) ^{1/4}$ $%
\forall m\in \left\{ 1,\cdot \cdot ,2N\right\} $. Therefore, by (\ref{3.10})
we have 
\begin{equation}
\underset{k=1,\cdot \cdot ,N}{\text{sup}}\left( \left\Vert \sigma
_{k}^{0}\right\Vert _{L^{2}\left( 0,T\right) }^{2}+\left\Vert \sigma
_{k}^{1}\right\Vert _{L^{2}\left( 0,T\right) }^{2}\right) \leq 2ce^{2\sqrt{%
\lambda _{N}}T}e^{\frac{c\sqrt{2}}{\varsigma }\sqrt{\lambda _{N}}}B\left(
T,\varsigma \left( \lambda _{N}\right) ^{-1/4}\right) \text{ .}  \tag{3.11}
\label{3.11}
\end{equation}

\bigskip

Our control given by $g\left( x,t\right) :=\sum\limits_{j=1,\cdot \cdot ,N}%
\left[ \alpha _{j}\sigma _{j}^{0}\left( t\right) +\beta _{j}\sigma
_{j}^{1}\left( t\right) \right] \Phi _{j}\left( x\right) $ where%
\begin{equation*}
\alpha _{j}:=\frac{\int_{0}^{1}\left( -u_{1}+\sqrt{\lambda _{j}}u_{0}\right)
\Phi _{j}}{\int_{\omega }\left\vert \Phi _{j}\right\vert ^{2}}=\frac{-b_{j}+%
\sqrt{\lambda _{j}}a_{j}}{\int_{\omega }\left\vert \Phi _{j}\right\vert ^{2}}%
\text{ and }\beta _{j}:=\frac{\int_{0}^{1}\left( -u_{1}-\sqrt{\lambda _{j}}%
u_{0}\right) \Phi _{j}}{\int_{\omega }\left\vert \Phi _{j}\right\vert ^{2}}=%
\frac{-b_{j}-\sqrt{\lambda _{j}}a_{j}}{\int_{\omega }\left\vert \Phi
_{j}\right\vert ^{2}}\text{ ,}
\end{equation*}%
satisfies%
\begin{equation}
\sum\limits_{j=1,\cdot \cdot ,N}\left( \alpha _{j}^{2}+\beta _{j}^{2}\right)
=2\sum\limits_{j=1,\cdot \cdot ,N}\frac{\left( \lambda
_{j}a_{j}^{2}+b_{j}^{2}\right) }{\left( \int_{\omega }\left\vert \Phi
_{j}\right\vert ^{2}\right) ^{2}}\leq \frac{2\left( 1+\lambda _{N}\right) }{%
\left( \underset{j=1,\cdot \cdot ,N}{\text{inf}}\int_{\omega }\left\vert
\Phi _{j}\right\vert ^{2}\right) ^{2}}\sum\limits_{j=1,\cdot \cdot ,N}\left(
a_{j}^{2}+b_{j}^{2}\right) \text{ .}  \tag{3.12}  \label{3.12}
\end{equation}

\bigskip

Combining the above estimates (\ref{3.11}) and (\ref{3.12}), there is some
universal constant $c$ independent of $T$ such that for any $N\geq 1$%
\begin{equation}
\begin{array}{ll}
\left\Vert g\right\Vert _{L^{2}\left( \left( 0,1\right) \times \left(
0,T\right) \right) }^{2} & =\displaystyle\sum\limits_{j=1,\cdot \cdot ,N}%
\displaystyle\int_{0}^{T}\left\vert \alpha _{j}\sigma _{j}^{0}\left(
t\right) +\beta _{j}\sigma _{j}^{1}\left( t\right) \right\vert ^{2}dt \\ 
& \leq 2\displaystyle\sum\limits_{j=1,\cdot \cdot ,N}\left( \alpha
_{j}^{2}+\beta _{j}^{2}\right) \underset{k=1,\cdot \cdot ,N}{\text{sup}}%
\left( \left\Vert \sigma _{k}^{0}\right\Vert _{L^{2}\left( 0,T\right)
}^{2}+\left\Vert \sigma _{k}^{1}\right\Vert _{L^{2}\left( 0,T\right)
}^{2}\right)  \\ 
& \leq \frac{8\left( 1+\lambda _{N}\right) }{\left( \underset{j=1,\cdot
\cdot ,N}{\text{inf}}\displaystyle\int_{\omega }\left\vert \Phi
_{j}\right\vert ^{2}\right) ^{2}}ce^{2\sqrt{\lambda _{N}}T}e^{\frac{c\sqrt{2}%
}{\varsigma }\sqrt{\lambda _{N}}}B\left( T,\varsigma \left( \lambda
_{N}\right) ^{-1/4}\right) \displaystyle\sum\limits_{j=1,\cdot \cdot
,N}\left( a_{j}^{2}+b_{j}^{2}\right) \text{ .}%
\end{array}
\tag{3.13}  \label{3.13}
\end{equation}

\bigskip

Recall that the bound%
\begin{equation*}
\left\Vert g\right\Vert _{L^{2}\left( \left( 0,1\right) \times \left(
0,T\right) \right) }^{2}:=\sum\limits_{j=1,\cdot \cdot
,N}\int_{0}^{T}\left\vert \alpha _{j}\sigma _{j}^{0}\left( t\right) +\beta
_{j}\sigma _{j}^{1}\left( t\right) \right\vert ^{2}dt\leq K\left\Vert \left(
u_{0},u_{1}\right) \right\Vert _{\left( L^{2}\left( 0,1\right) \right)
^{2}}^{2}:=K\sum\limits_{j=1,\cdot \cdot ,N}\left( a_{j}^{2}+b_{j}^{2}\right)
\end{equation*}%
will imply that the solution $\varphi $ of (\ref{2.1}) satisfies 
\begin{equation*}
\left\Vert \varphi \left( \cdot ,T\right) \right\Vert _{L^{2}\left(
0,1\right) }^{2}+\left\Vert \partial _{t}\varphi \left( \cdot ,T\right)
\right\Vert _{L^{2}\left( 0,1\right) }^{2}\leq K\int_{0}^{T}\int_{\omega
}\left\vert \varphi \right\vert ^{2}\text{ .}
\end{equation*}

\bigskip

Now our aim is to bound the quantity 
\begin{equation*}
\frac{8\left( 1+\lambda _{N}\right) }{\left( \underset{j=1,\cdot \cdot ,N}{%
\text{inf}}\int_{\omega }\left\vert \Phi _{j}\right\vert ^{2}\right) ^{2}}%
ce^{2\sqrt{\lambda _{N}}T}e^{\frac{c\sqrt{2}}{\varsigma }\sqrt{\lambda _{N}}%
}B\left( T,\varsigma \left( \lambda _{N}\right) ^{-1/4}\right)
\end{equation*}%
appearing in (\ref{3.13}) in order to get the cost $K$.

\bigskip

First, by Proposition 2.1, $\frac{1}{\left( \underset{j=1,\cdot \cdot ,N}{%
\text{inf}}\int_{\omega }\left\vert \Phi _{j}\right\vert ^{2}\right) ^{2}}%
\leq \frac{1}{\rho ^{2}(2-\alpha )^{2}}$. Next, recall that $\varsigma =$min$%
\left( \frac{\gamma (2-\alpha )}{2\sqrt{2}},\frac{2\sqrt{\lambda _{1}}}{1+%
\sqrt{2}}\right) $ and since $\alpha \in \left[ 0,2\right) $ with (\ref{1.2}%
), we have that $\underline{c}\gamma (2-\alpha )\leq \varsigma \leq \frac{1}{%
\underline{c}}$ where $\underline{c}$ is a positive constant independent on $%
\alpha \in \lbrack 0,2)$. Finally, the estimate of $B\left( T,r\right) $ in
Theorem 3.1 
\begin{equation*}
B\left( T,r\right) =\left\{ 
\begin{array}{cc}
\left( \frac{1}{T}+\frac{1}{T^{2}r^{2}}\right) e^{c\frac{1}{Tr^{2}}} & \text{%
if }T\leq \frac{1}{r^{2}} \\ 
cr^{2} & \text{if }T\geq \frac{1}{r^{2}}%
\end{array}%
\right. \leq \left\{ 
\begin{array}{cc}
\left( 1+\frac{1}{c}\right) \frac{1}{T}e^{2c\frac{1}{Tr^{2}}} & \text{if }%
T\leq \frac{1}{r^{2}} \\ 
cr^{2} & \text{if }T\geq \frac{1}{r^{2}}%
\end{array}%
\right. \leq \left( (1+\frac{1}{c})\frac{1}{T}+cr^{2}\right) e^{2c\frac{1}{%
Tr^{2}}}
\end{equation*}%
leads to the bound 
\begin{eqnarray*}
B(T,\varsigma \left( \lambda _{N}\right) ^{-1/4}) &\leq &\left( (1+\frac{1}{c%
})\frac{1}{T}+c\frac{\varsigma ^{2}}{\sqrt{\lambda _{N}}}\right) e^{2c\frac{%
\sqrt{\lambda _{N}}}{T\varsigma ^{2}}} \\
&\leq &C(1+\frac{1}{T})e^{C\frac{\sqrt{\lambda _{N}}}{T(2-\alpha )^{2}}}
\end{eqnarray*}%
for some $C>0$ independent on $N>0$, $\alpha \in \lbrack 0,2)$ and $T>0$.
Therefore, by (\ref{3.13}) one can conclude that 
\begin{equation*}
\left\Vert g\right\Vert _{L^{2}\left( \left( 0,1\right) \times \left(
0,T\right) \right) }^{2}\leq \frac{C\left( 1+\lambda _{N}\right) }{\rho
^{2}(2-\alpha )^{2}}e^{C\sqrt{\lambda _{N}}T}e^{C\frac{\sqrt{\lambda _{N}}}{%
\gamma \left( 2-\alpha \right) }}C(1+\frac{1}{T})e^{C\frac{\sqrt{\lambda _{N}%
}}{T\gamma ^{2}(2-\alpha )^{2}}}\sum\limits_{j=1,\cdot \cdot ,N}\left(
a_{j}^{2}+b_{j}^{2}\right) \text{ ,}
\end{equation*}%
which gives, using $\frac{1}{\gamma \left( 2-\alpha \right) }\leq T+\frac{1}{%
T\gamma ^{2}(2-\alpha )^{2}}$ that 
\begin{equation*}
\left\Vert g\right\Vert _{L^{2}\left( \left( 0,1\right) \times \left(
0,T\right) \right) }^{2}\leq \frac{C\left( 1+\lambda _{N}\right) }{\rho
^{2}(2-\alpha )^{2}}\left( 1+\frac{1}{T}\right) e^{C\sqrt{\lambda _{N}}%
\left( T+\frac{1}{T\gamma ^{2}(2-\alpha )^{2}}\right)
}\sum\limits_{j=1,\cdot \cdot ,N}\left( a_{j}^{2}+b_{j}^{2}\right) \text{ .}
\end{equation*}%
By the cost estimate in Lemma 3.1, we obtain that for any $\varphi $
solution of (\ref{2.1}) and any $N\geq 1$ 
\begin{equation*}
\left\Vert \varphi \left( \cdot ,T\right) \right\Vert _{L^{2}\left(
0,1\right) }^{2}+\left\Vert \partial _{t}\varphi \left( \cdot ,T\right)
\right\Vert _{L^{2}\left( 0,1\right) }^{2}\leq \frac{C\left( 1+\lambda
_{N}\right) }{\rho ^{2}(2-\alpha )^{2}}\left( 1+\frac{1}{T}\right) e^{C\sqrt{%
\lambda _{N}}\left( T+\frac{1}{T\gamma ^{2}(2-\alpha )^{2}}\right)
}\int_{0}^{T}\int_{\omega }\left\vert \varphi \right\vert ^{2}\text{ ,}
\end{equation*}%
where $C>0$ does not depend on $\left( N,T,\alpha \right) $. This completes
the proof of Proposition 2.3.

\bigskip

%Finally, we choose $T=\frac{\sqrt{%
%\lambda _{1}}}{\text{cst}^{2}}$ in order that $T\leq \frac{1}{r^{2}}:=\frac{%
%\sqrt{\lambda _{N}}}{\text{cst}^{2}}$ and

\bigskip

\section{Elliptic observation by Carleman techniques (proof of Proposition
2.2)}

\bigskip

In this section, we shall prove Proposition 2.2. Let $0<a<b<1$ and $\Omega
=\left( a,b\right) \times \left( 0,T\right) $. We set $\left( x,t\right)
=\left( x_{1},x_{2}\right) \in \Omega $, and for $\alpha \in \lbrack 0,2)$,
introduce 
\begin{equation*}
Q=-\partial _{t}^{2}-\mathcal{P}=-\nabla \cdot (A(x_{1},x_{2})\nabla \cdot
),\quad A(x_{1},x_{2})=%
\begin{pmatrix}
x_{1}^{\alpha } & 0 \\ 
0 & 1%
\end{pmatrix}%
,\quad \nabla =%
\begin{pmatrix}
\partial _{1} \\ 
\partial _{2}%
\end{pmatrix}%
\text{ .}
\end{equation*}%
Note that there exists $C_{0}>0$ such that 
\begin{equation}
\left\Vert A\right\Vert _{W^{3,\infty }(\Omega )}\leq C_{0},\quad
A(x_{1},x_{2})\xi \cdot \xi \geq \frac{1}{C_{0}}|\xi |^{2},\quad \forall \xi
\in \mathbb{R}^{2},\forall (x_{1},x_{2})\in \Omega \text{ ,}  \tag{4.1}
\label{4.1}
\end{equation}%
where $C_{0}>0$ is independent on $\alpha \in \lbrack 0,2)$. We set 
\begin{equation*}
v=e^{\tau \phi }\chi z
\end{equation*}%
where $\tau >0$, $z\in H^{2}\left( \Omega \right) $, $\chi \left(
x_{1},x_{2}\right) =\chi _{1}\left( x_{1}\right) \chi _{2}\left(
x_{2}\right) $ with%
\begin{equation*}
\left\{ 
\begin{array}{lll}
\chi _{1}\in C_{0}^{\infty }\left( a,b\right) \text{,} & 0\leq \chi _{1}\leq
1\text{,} & \chi _{1}=1\text{ on }\left( \frac{3a+b}{4},\frac{a+3b}{4}\right)
\\ 
\chi _{2}\in C^{\infty }\left( 0,T\right) \text{,} & 0\leq \chi _{2}\leq 1%
\text{,} & \chi _{2}=1\text{ on }\left( 0,\frac{T}{3}\right) \text{ and }%
\chi _{2}=0\text{ on }\left( \frac{2T}{3},T\right)%
\end{array}%
\right.
\end{equation*}%
and we shall consider weight functions $\phi \in C^{\infty }(\overline{%
\Omega })$ of the form 
\begin{equation}
\phi (x_{1},x_{2})=e^{\lambda \psi (x_{1},x_{2})},\quad \lambda >0,\quad
\psi \in C^{\infty }(\overline{\Omega }),\quad \nabla \psi \neq 0\text{ on }%
\overline{\Omega }\text{ .}  \tag{4.2}  \label{4.2}
\end{equation}%
Here, we give explicitely $\psi $ as follows%
\begin{equation}
\psi (x_{1},x_{2})=-\left( x_{1}-x_{0}\right) ^{2k}-\beta ^{2k}\left(
x_{2}+1\right) ^{2k}  \tag{4.3}  \label{4.3}
\end{equation}%
where $x_{0}=\frac{a+b}{2}$, $\beta =\frac{2}{3}\left( \frac{b-a}{T+4}%
\right) $ and $k=$max$\left( \text{ln}2\text{/ln}\left( \left( 4T+12\right)
/\left( 3T+12\right) \right) ;\text{ln}2\text{/ln}\left( 3/2\right) \right) $%
.

\bigskip

We set 
\begin{equation*}
Q_{\phi }=e^{\tau \phi }Qe^{-\tau \phi }\text{ .}
\end{equation*}%
We have $Q_{\phi }v=\mathcal{S}v+\mathcal{A}v+\mathcal{R}v$ with 
\begin{equation*}
\mathcal{S}v=-\nabla \cdot (A\nabla v)-\tau ^{2}A\nabla \phi \cdot \nabla
\phi v,\quad \mathcal{A}v=2\tau A\nabla \phi \cdot \nabla v+2\tau \nabla
\cdot (A\nabla \phi )v,\quad \mathcal{R}v=-\tau \nabla \cdot (A\nabla \phi )v%
\text{ ,}
\end{equation*}%
which gives $\left\Vert Q_{\phi }v-\mathcal{R}v\right\Vert _{L^{2}(\Omega
)}^{2}=\left\Vert \mathcal{S}v\right\Vert _{L^{2}(\Omega )}^{2}+\left\Vert 
\mathcal{A}v\right\Vert _{L^{2}(\Omega )}^{2}+2(\mathcal{S}v,\mathcal{A}%
v)_{L^{2}(\Omega )}$. Note that $0\leq \left\Vert Q_{\phi }v-\mathcal{R}%
v\right\Vert _{L^{2}(\Omega )}^{2}-2(\mathcal{S}v,\mathcal{A}%
v)_{L^{2}(\Omega )}$ implies 
\begin{equation}
\left( \mathcal{S}v,\mathcal{A}v\right) _{L^{2}(\Omega )}\leq \left\Vert
Q_{\phi }v\right\Vert _{L^{2}(\Omega )}^{2}+\left\Vert \mathcal{R}%
v\right\Vert _{L^{2}(\Omega )}^{2}\text{ .}  \tag{4.4}  \label{4.4}
\end{equation}

\bigskip

Now we compute $\left( \mathcal{S}v,\mathcal{A}v\right) _{L^{2}(\Omega )}$:
By integration by parts, one has with standard summation notations and $%
A=\left( A_{ij}\right) _{1\leq i,j\leq 2}$,%
\begin{equation*}
\begin{array}{ll}
\left( \mathcal{S}v,\mathcal{A}v\right) _{L^{2}(\Omega )} & =2\tau %
\displaystyle\int_{\Omega }A\nabla ^{2}vA\nabla \phi \cdot \nabla v+2\tau %
\displaystyle\int_{\Omega }A\nabla ^{2}\phi A\nabla v\cdot \nabla v \\ 
& \text{\quad }+2\tau \displaystyle\int_{\Omega }A_{ij}\partial
_{x_{i}}v\partial _{x_{\ell }}v\partial _{x_{j}}A_{k\ell }\partial
_{x_{k}}\phi \\ 
& \text{\quad }+2\tau \displaystyle\int_{\Omega }\left( A\nabla v\cdot
\nabla v\right) \nabla \cdot \left( A\nabla \phi \right) +2\tau \displaystyle%
\int_{\Omega }A\nabla v\cdot \nabla \left( \nabla \cdot \left( A\nabla \phi
\right) \right) v \\ 
& \text{\quad }+\tau ^{3}\displaystyle\int_{\Omega }\left[ A\nabla \phi
\cdot \nabla \left( A\nabla \phi \cdot \nabla \phi \right) -\left( A\nabla
\phi \cdot \nabla \phi \right) \left( \nabla \cdot \left( A\nabla \phi
\right) \right) \right] \left\vert v\right\vert ^{2} \\ 
& \text{\quad }+2\tau \displaystyle\int_{\partial \Omega }\left( A\nabla
v\cdot n\right) \left( A\nabla \phi \cdot \nabla v+\left( \nabla \cdot
\left( A\nabla \phi \right) \right) v\right) -\tau ^{3}\displaystyle%
\int_{\partial \Omega }\left( A\nabla \phi \cdot \nabla \phi \right) \left(
A\nabla \phi \cdot n\right) \left\vert v\right\vert ^{2}\text{ .}%
\end{array}%
\end{equation*}%
But by one integration by parts%
\begin{equation*}
\begin{array}{ll}
\displaystyle\int_{\Omega }A\nabla ^{2}vA\nabla \phi \cdot \nabla v & =%
\displaystyle\frac{1}{2}\int_{\partial \Omega }\left( A\nabla v\cdot \nabla
v\right) \left( A\nabla \phi \cdot n\right) -\displaystyle\frac{1}{2}%
\int_{\Omega }\partial _{x_{\ell }}A_{ij}\partial _{x_{j}}v\partial
_{x_{i}}vA_{k\ell }\partial _{x_{k}}\phi \\ 
& \quad -\displaystyle\frac{1}{2}\int_{\Omega }\left( A\nabla v\cdot \nabla
v\right) \nabla \cdot \left( A\nabla \phi \right) \text{ .}%
\end{array}%
\end{equation*}%
Therefore, 
\begin{equation*}
\begin{array}{ll}
\left( \mathcal{S}v,\mathcal{A}v\right) _{L^{2}(\Omega )} & =2\tau %
\displaystyle\int_{\Omega }A\nabla ^{2}\phi A\nabla v\cdot \nabla v+\tau %
\displaystyle\int_{\Omega }\left( A\nabla v\cdot \nabla v\right) \nabla
\cdot \left( A\nabla \phi \right) \\ 
& \text{\quad }+\tau ^{3}\displaystyle\int_{\Omega }\left[ A\nabla \phi
\cdot \nabla \left( A\nabla \phi \cdot \nabla \phi \right) -\left( A\nabla
\phi \cdot \nabla \phi \right) \left( \nabla \cdot \left( A\nabla \phi
\right) \right) \right] \left\vert v\right\vert ^{2} \\ 
& \text{\quad }+R_{1}+R_{2}%
\end{array}%
\end{equation*}%
with%
\begin{equation*}
\begin{array}{ll}
R_{1} & =2\tau \displaystyle\int_{\Omega }A_{ij}\partial _{x_{i}}v\partial
_{x_{\ell }}v\partial _{x_{j}}A_{k\ell }\partial _{x_{k}}\phi -\tau %
\displaystyle\int_{\Omega }\partial _{x_{\ell }}A_{ij}\partial
_{x_{j}}v\partial _{x_{i}}vA_{k\ell }\partial _{x_{k}}\phi \\ 
& \text{\quad }+2\tau \displaystyle\int_{\Omega }A\nabla v\cdot \nabla
\left( \nabla \cdot \left( A\nabla \phi \right) \right) v\text{ ,}%
\end{array}%
\end{equation*}%
\begin{equation*}
\begin{array}{ll}
R_{2} & =-2\tau \displaystyle\int_{\partial \Omega }\left( A\nabla v\cdot
n\right) \left( A\nabla \phi \cdot \nabla v\right) +\tau \displaystyle%
\int_{\partial \Omega }\left( A\nabla v\cdot \nabla v\right) \left( A\nabla
\phi \cdot n\right) \\ 
& \text{\quad }-2\tau \displaystyle\int_{\partial \Omega }\left( A\nabla
v\cdot n\right) \left( \nabla \cdot \left( A\nabla \phi \right) \right)
v-\tau ^{3}\displaystyle\int_{\partial \Omega }\left( A\nabla \phi \cdot
\nabla \phi \right) \left( A\nabla \phi \cdot n\right) \left\vert
v\right\vert ^{2}\text{ ,}%
\end{array}%
\end{equation*}%
where $n$ is the outward normal vector to $\partial \Omega $.

\bigskip

Notice that from the form of $A$ and $\phi $ given by (\ref{4.1}) and (\ref%
{4.2}), we have the existence of $C_{1}>0$ independent on $\alpha \in
\lbrack 0,2)$ such that for $\tau >0$ sufficiently large 
\begin{equation*}
\left\vert R_{1}\right\vert \leq C_{1}\left( (\tau ^{1/2}\lambda ^{2}+\tau
\lambda )\left\Vert \phi ^{1/2}\nabla v\right\Vert _{L^{2}(\Omega
)}^{2}+\tau ^{3/2}\lambda ^{4}\left\Vert \phi ^{1/2}v\right\Vert
_{L^{2}(\Omega )}^{2}\right) \text{ .}
\end{equation*}%
Note also that from the form of $A$ and $\phi $ given by (\ref{4.1}) and (%
\ref{4.2}), we have 
\begin{equation*}
A\nabla ^{2}\phi A\nabla v\cdot \nabla v=\lambda ^{2}\phi (A\nabla \psi
\cdot \nabla v)^{2}+\lambda \phi A\nabla ^{2}\psi A\nabla v\cdot \nabla
v\geq -C_{2}\lambda \phi |\nabla v|^{2}\text{ ,}
\end{equation*}%
and 
\begin{align*}
\tau \int_{\Omega }(A\nabla v\cdot \nabla v)\nabla \cdot (A\nabla \phi )&
=\tau \int_{\Omega }(A\nabla v\cdot \nabla v)\phi (\lambda \nabla \cdot
(A\nabla \psi )+\lambda ^{2}A\nabla \psi \cdot \nabla \psi )) \\
& \geq C_{2}\tau \lambda ^{2}\left\Vert \phi ^{1/2}\nabla v\right\Vert
_{L^{2}(\Omega )}^{2}-C_{3}\tau \lambda \left\Vert \phi ^{1/2}\nabla
v\right\Vert _{L^{2}(\Omega )}^{2} \\
& \geq C_{4}\tau \lambda ^{2}\left\Vert \phi ^{1/2}\nabla v\right\Vert
_{L^{2}(\Omega )}^{2}\text{ ,}
\end{align*}%
for $\lambda >0$ chosen sufficiently large (independently on $\alpha \in
\lbrack 0,2)$, and where the constants $C_{2},C_{3},C_{4}>0$ are independent
on $\alpha \in \lbrack 0,2)$. Arguing in the same way, there exist constants 
$C_{5}>0$ and $\lambda _{0}>0$ such that for all $\alpha \in \lbrack 0,2)$
and for all $\lambda >\lambda _{0}$, 
\begin{equation*}
\tau ^{3}\int_{\Omega }\left[ A\nabla \phi \cdot \nabla \left( A\nabla \phi
\cdot \nabla \phi \right) -\left( A\nabla \phi \cdot \nabla \phi \right)
\left( \nabla \cdot \left( A\nabla \phi \right) \right) \right] \left\vert
v\right\vert ^{2}\geq C_{5}\tau ^{3}\lambda ^{4}\left\Vert \phi
^{3/2}v\right\Vert _{L^{2}(\Omega )}^{2}\text{ .}
\end{equation*}%
Summing up, (\ref{4.4}) becomes 
\begin{multline*}
C_{5}\tau ^{3}\lambda ^{4}\left\Vert \phi ^{3/2}v\right\Vert _{L^{2}(\Omega
)}^{2}+C_{4}\tau \lambda ^{2}\left\Vert \phi ^{1/2}\nabla v\right\Vert
_{L^{2}(\Omega )}^{2}+R_{2} \\
\leq C_{1}\left( (\tau ^{1/2}\lambda ^{2}+\tau \lambda )\left\Vert \phi
^{1/2}\nabla v\right\Vert _{L^{2}(\Omega )}^{2}+\tau ^{3/2}\lambda
^{4}\left\Vert \phi ^{1/2}v\right\Vert _{L^{2}(\Omega )}^{2}\right)
+\left\Vert Q_{\phi }v\right\Vert _{L^{2}(\Omega )}^{2}+\left\Vert \mathcal{R%
}v\right\Vert _{L^{2}(\Omega )}^{2}\text{ ,}
\end{multline*}%
where the constants are independent on $\alpha \in \lbrack 0,2)$. Fixing $%
\lambda >\lambda _{0}$ large, and then taking $\tau >\tau _{0}$ sufficiently
large (constants may depend on $\lambda $ from now), we obtain the existence
of $C_{6}>0$ such that 
\begin{equation*}
C_{6}\tau ^{3}\left\Vert v\right\Vert _{L^{2}(\Omega )}^{2}+C_{6}\tau
\left\Vert \nabla v\right\Vert _{L^{2}(\Omega )}^{2}+R_{2}\leq \left\Vert
Q_{\phi }v\right\Vert _{L^{2}(\Omega )}^{2}+\left\Vert \mathcal{R}%
v\right\Vert _{L^{2}(\Omega )}^{2}\text{ .}
\end{equation*}%
Next, one can see that from the form of $A$ and $\phi $, there is $C_{7}>0$
such that for all $\alpha \in \lbrack 0,2)$, 
\begin{equation*}
\left\Vert \mathcal{R}v\right\Vert _{L^{2}(\Omega )}^{2}\leq C_{7}\tau
^{2}\left\Vert v\right\Vert _{L^{2}(\Omega )}^{2}\text{ .}
\end{equation*}%
Therefore, taking $\tau >0$ sufficiently large yields the existence of $%
C_{8}>0$ such that 
\begin{equation}
C_{8}\left( \tau ^{3}\left\Vert v\right\Vert _{L^{2}(\Omega )}^{2}+\tau
\left\Vert \nabla v\right\Vert _{L^{2}(\Omega )}^{2}\right) +R_{2}\leq
\left\Vert Q_{\phi }v\right\Vert _{L^{2}(\Omega )}^{2}\text{ .}  \tag{4.5}
\label{4.5}
\end{equation}

\bigskip

Now we treat the boundary term $R_{2}$: Since $v=A\nabla v\cdot n=0$ on $%
\partial \Omega \left\backslash \Gamma \right. $ where $\Gamma =\left\{
\left( x_{1},0\right) ;x_{1}\in \left( a,b\right) \right\} $, one can deduce
that 
\begin{equation*}
\begin{array}{ll}
R_{2} & =\tau \displaystyle\int_{a}^{b}\partial _{2}\phi \left\vert \partial
_{2}v\left( x_{1},0\right) \right\vert ^{2}dx_{1} \\ 
& \text{\quad }+2\tau \displaystyle\int_{a}^{b}x_{1}^{\alpha }\partial
_{1}\phi \partial _{1}v\left( x_{1},0\right) \partial _{2}v\left(
x_{1},0\right) dx_{1}-\tau \displaystyle\int_{a}^{b}x_{1}^{\alpha }\partial
_{2}\phi \left\vert \partial _{1}v\left( x_{1},0\right) \right\vert
^{2}dx_{1} \\ 
& \text{\quad }+2\tau \displaystyle\int_{a}^{b}\left( \nabla \cdot \left(
A\nabla \phi \right) \right) v\left( x_{1},0\right) \partial _{2}v\left(
x_{1},0\right) dx_{1}+\tau ^{3}\displaystyle\int_{a}^{b}\left( A\nabla \phi
\cdot \nabla \phi \right) \partial _{2}\phi \left\vert v\left(
x_{1},0\right) \right\vert ^{2}dx_{1}\text{ .}%
\end{array}%
\end{equation*}%
which gives the existence of $C_{9}>0$ independent on $\alpha \in \lbrack
0,2)$ such that for any $\tau >0$ sufficiently large 
\begin{equation*}
\left\vert R_{2}\right\vert \leq C_{9}\left( \tau \left\Vert \partial
_{2}v\left( \cdot ,0\right) \right\Vert _{L^{2}(a,b)}^{2}+\tau
^{3}\left\Vert v\left( \cdot ,0\right) \right\Vert _{H^{1}(a,b)}^{2}\right) 
\text{ .}
\end{equation*}

\bigskip

Finally, by (\ref{4.5}) we have for any $\tau >\tau _{0}$ with $\tau _{0}>1$%
, the following inequality%
\begin{equation}
C_{8}\left( \tau ^{3}\left\Vert v\right\Vert _{L^{2}(\Omega )}^{2}+\tau
\left\Vert \nabla v\right\Vert _{L^{2}(\Omega )}^{2}\right) \leq \left\Vert
Q_{\phi }v\right\Vert _{L^{2}(\Omega )}^{2}+C_{9}\left( \tau \left\Vert
\partial _{2}v\left( \cdot ,0\right) \right\Vert _{L^{2}(a,b)}^{2}+\tau
^{3}\left\Vert v\left( \cdot ,0\right) \right\Vert _{H^{1}(a,b)}^{2}\right) 
\text{ .}  \tag{4.6}  \label{4.6}
\end{equation}

\bigskip

Let $U=\left( \frac{2a+b}{3},\frac{a+2b}{3}\right) \times \left( 0,\frac{T}{4%
}\right) $, $W_{1}=\left( \left[ a,\frac{3a+b}{4}\right] \cup \left[ \frac{%
a+3b}{4},b\right] \right) \times \left[ 0,\frac{2T}{3}\right] $, $W_{2}=%
\left[ a,b\right] \times \left[ \frac{T}{3},\frac{2T}{3}\right] $ and $%
W=W_{1}\cup W_{2}$. We have supp$\nabla \chi =W$ and $\chi =1$ in $U$.

\bigskip

Coming back to the function $z$ where $v=e^{\tau \phi }\chi z$, $Q_{\phi
}v=e^{\tau \phi }Q\left( \chi z\right) =e^{\tau \phi }\left( \chi Qz+\left[
Q,\chi \right] z\right) $ where the bracket $[Q,\chi ]=-\partial
_{t}^{2}\chi -2(\partial _{t}\chi )\partial _{t}-x^{\alpha }(\partial
_{x}^{2}\chi )-2(\partial _{x}\chi )x^{\alpha }\partial _{x}-\alpha
(\partial _{x}\chi )x^{\alpha -1}$ is a differential operator of order one,
supported in $W$, which is away from a neighborhood of the degeneracy $%
\{x=0\}$. From (\ref{4.6}) and taking any $\tau $ sufficiently large yields%
\begin{eqnarray*}
\tau ^{3}\left\Vert e^{\tau \phi }z\right\Vert _{L^{2}(U)}^{2}+\tau
\left\Vert e^{\tau \phi }\nabla z\right\Vert _{L^{2}(U)}^{2} &\leq &C\left(
\left\Vert e^{\tau \phi }\chi Qz\right\Vert _{L^{2}(\Omega )}^{2}+\tau
\left\Vert e^{\tau \phi }z\right\Vert _{L^{2}(W)}^{2}+\left\Vert e^{\tau
\phi }\nabla z\right\Vert _{L^{2}(W)}^{2}\right) \\
&&+C\left( \tau \left\Vert e^{\tau \phi \left( \cdot ,0\right) }\partial
_{t}z\left( \cdot ,0\right) \right\Vert _{L^{2}(a,b)}^{2}+\tau
^{5}\left\Vert e^{\tau \phi \left( \cdot ,0\right) }z\left( \cdot ,0\right)
\right\Vert _{H^{1}(a,b)}^{2}\right) \text{ .}
\end{eqnarray*}%
Let $D=\underset{\Omega }{\text{max}}\phi $, $D_{W}=\underset{W}{\text{max}}%
\phi $, $D_{0}=\underset{\left( a,b\right) }{\text{max}}\phi \left( \cdot
,0\right) $ and $D_{U}=\underset{U}{\text{min}}\phi $. We have for any $\tau
>\tau _{0}$ sufficiently large 
\begin{eqnarray*}
e^{2\tau D_{U}}\left( \left\Vert z\right\Vert _{L^{2}(U)}^{2}+\left\Vert
\nabla z\right\Vert _{L^{2}(U)}^{2}\right) &\leq &Ce^{2\tau D}\left\Vert
Qz\right\Vert _{L^{2}(\Omega )}^{2}+Ce^{2\tau D_{K}}\left( \tau \left\Vert
z\right\Vert _{L^{2}(W)}^{2}+\left\Vert \nabla z\right\Vert
_{L^{2}(W)}^{2}\right) \\
&&+Ce^{2\tau D_{0}}\left( \tau \left\Vert \partial _{t}z\left( \cdot
,0\right) \right\Vert _{L^{2}(a,b)}^{2}+\tau ^{5}\left\Vert z\left( \cdot
,0\right) \right\Vert _{H^{1}(a,b)}^{2}\right) \text{ .}
\end{eqnarray*}%
Our choice of $\psi $ given by (\ref{4.3}) allows to get $D>D_{U}$ and $%
D_{0}>D_{U}>D_{K}$. Indeed, by a straightforward computation, 
\begin{equation*}
\left\{ 
\begin{array}{rl}
\underset{W_{1}}{\text{max}}\psi -\underset{U}{\text{min}}\psi & \leq
-\left( \frac{b-a}{4}\right) ^{2k}-\beta ^{2k}+\left( \frac{b-a}{6}\right)
^{2k}+\beta ^{2k}\left( \frac{T}{4}+1\right) ^{2k} \\ 
& =\beta ^{2k}\left( -1+\left( \frac{3}{8}\left( T+4\right) \right)
^{2k}\left( -1+2\left( \frac{2}{3}\right) ^{2k}\right) \right) <0\text{ ,}
\\ 
\underset{W_{2}}{\text{max}}\psi -\underset{U}{\text{min}}\psi & \leq -\beta
^{2k}\left( \frac{T}{3}+1\right) ^{2k}+\left( \frac{b-a}{6}\right)
^{2k}+\beta ^{2k}\left( \frac{T}{4}+1\right) ^{2k} \\ 
& =\beta ^{2k}\left( \frac{1}{3}\left( T+4\right) \right) ^{2k}\left(
-1+2\left( \frac{3T+12}{4T+12}\right) ^{2k}\right) <0\text{ .}%
\end{array}%
\right.
\end{equation*}

Using $W\subset \overline{\Omega }$ and optimizing with respect to $\tau $
yield the desired interpolation estimate (see e.g. \cite{R} or \cite[Lemma
5.4, page 189]{LRLeR1}). This completes the proof of%
\begin{equation*}
\left\Vert \varphi \right\Vert _{H^{1}\left( \left( \frac{2a+b}{3},\frac{a+2b%
}{3}\right) \times \left( 0,\frac{T}{4}\right) \right) }\leq c\left\Vert
\varphi \right\Vert _{H^{1}\left( \left( a.b\right) \times \left( 0,T\right)
\right) }^{1-\delta }\left( \left\Vert \varphi _{0}\right\Vert _{H^{1}\left(
a,b\right) }+\left\Vert \varphi _{1}\right\Vert _{L^{2}\left( a,b\right)
}\right) ^{\delta }\text{ ,}
\end{equation*}%
since $Q\varphi =0$.

\bigskip

\bigskip

\section{Observability estimate for the eigenfunctions (proof of Proposition
2.1)}

\bigskip

In this section we aim to prove Proposition 2.1. Given $0<a<b<1$, we will
use the notation $X\lesssim Y$ , or $Y\gtrsim X$ to denote the bound $%
\left\vert X\right\vert \leq cY$ for some constant $c>0$ only dependent on $%
\left( a,b\right) $.

\bigskip

Cannarsa, Martinez and Vancostenoble proved (see \cite{CMV4} proposition
2.15 at page 10) that 
\begin{equation*}
\forall \alpha \in \left[ 1,2\right) \text{ \quad }\forall j\geq 1\text{
\quad }||\Phi _{j}||_{L^{2}(a,b)}^{2}\gtrsim 2-\alpha \text{ .}
\end{equation*}%
In this section, we extend this result to $\alpha \in \left[ 0,2\right) $.
To this end, we focus on the case $\alpha \in \left[ 0,1\right) $ and apply
the following observability estimate.

\bigskip

Proposition 5.1 .-\textit{\ For all }$\sigma \in \mathbb{R}$\textit{, for
all }$\alpha \in \left[ 0,1\right) $\textit{, for all }$\vartheta \in D(%
\mathcal{P})$ 
\begin{equation*}
\sigma ^{2}\left\Vert \vartheta \right\Vert _{L^{2}(0,1)}^{2}+\left\Vert
x^{\alpha /2}\vartheta ^{\prime }\right\Vert _{L^{2}(0,1)}^{2}\lesssim
\left( \left\Vert \left( \mathcal{P}-\sigma ^{2}\right) \vartheta
\right\Vert _{L^{2}(0,1)}^{2}+(1+\sigma ^{2})\left\Vert \vartheta
\right\Vert _{L^{2}(a,b)}^{2}\right) \text{ .}
\end{equation*}

\bigskip

Since $\Phi _{j}\in D(\mathcal{P})$ is the normalized eigenfunctions of $%
\mathcal{P}$ associated with an eigenvalue $\lambda _{j}$, $j\in \mathbb{N}%
^{\ast }$. Applying Proposition 5.1 with $\vartheta =\Phi _{j}$ and $\sigma
^{2}=\lambda _{j}$, we obtain 
\begin{equation*}
\frac{\lambda _{j}}{1+\lambda _{j}}\lesssim \left\Vert \Phi _{j}\right\Vert
_{L^{2}(a,b)}^{2}\text{ .}
\end{equation*}%
Using $\frac{\lambda _{1}}{1+\lambda _{1}}\leq \frac{\lambda _{j}}{1+\lambda
_{j}}$ and (\ref{1.2}), one can deduce that%
\begin{equation*}
\forall \alpha \in \left[ 0,1\right) \text{ \quad }\forall j\geq 1\text{
\quad }\left\Vert \Phi _{j}\right\Vert _{L^{2}(a,b)}^{2}\gtrsim 1\geq \frac{1%
}{2}\left( 2-\alpha \right) \text{ .}
\end{equation*}%
This completes the proof of Proposition 2.1.

\bigskip

Now, we prove Proposition 5.1. Before proceeding to the proof we need two
lemmas.

\bigskip

Lemma 5.1 .- \textit{There exists }$C>0$\textit{\ such that for all }$\sigma
\in \mathbb{R}$\textit{, for all }$\alpha \in \lbrack 0,1)$\textit{,} 
\begin{equation*}
\sigma ^{2}\left\Vert \vartheta \right\Vert _{L^{2}(0,1)}^{2}+\left\Vert
x^{\alpha /2}\vartheta ^{\prime }\right\Vert _{L^{2}(0,1)}^{2}\leq C\left(
\left\Vert \left( \mathcal{P}-\sigma ^{2}\right) \vartheta \right\Vert
_{L^{2}(0,1)}^{2}+\left\vert \vartheta ^{\prime }(1)\right\vert ^{2}\right) 
\text{ ,}
\end{equation*}%
\textit{for all }$\vartheta \in D(\mathcal{P})$.

\bigskip

Lemma 5.2 - \textit{There exists }$C>0$\textit{\ such that for all }$\sigma
\in \mathbb{R}$\textit{, for all }$\alpha \in \lbrack 0,1)$\textit{,} 
\begin{equation*}
\left\vert \vartheta ^{\prime }(1)\right\vert ^{2}\leq C\left( \left\Vert
\left( \mathcal{P}-\sigma ^{2}\right) \vartheta \right\Vert _{L^{2}\left(
0,1\right) }^{2}+\left\Vert \vartheta \right\Vert _{H^{1}\left( \frac{3a+b}{4%
},\frac{a+3b}{4}\right) }^{2}\right) \text{ ,}
\end{equation*}%
\textit{for all }$\vartheta \in D(\mathcal{P})$.

\bigskip

Proof of Proposition 5.1 .- By Lemmata 5.1 and 5.2, 
\begin{equation*}
\sigma ^{2}\left\Vert \vartheta \right\Vert _{L^{2}(0,1)}^{2}+\left\Vert
x^{\alpha /2}\vartheta ^{\prime }\right\Vert _{L^{2}(0,1)}^{2}\lesssim
C\left( \left\Vert \left( \mathcal{P}-\sigma ^{2}\right) \vartheta
\right\Vert _{L^{2}\left( 0,1\right) }^{2}+\left\Vert \vartheta \right\Vert
_{H^{1}\left( \frac{3a+b}{4},\frac{a+3b}{4}\right) }^{2}\right) \text{ .}
\end{equation*}%
Let $\chi \in C_{0}^{\infty }\left( 0,1\right) $ such that $0\leq \chi \leq
1 $ and $\chi =1$ on $[\frac{3a+b}{4},\frac{a+3b}{4}]$. We have 
\begin{align*}
\left\Vert \vartheta \right\Vert _{H^{1}\left( \frac{3a+b}{4},\frac{a+3b}{4}%
\right) }^{2}& =\left\Vert \chi \vartheta \right\Vert _{H^{1}\left( \frac{%
3a+b}{4},\frac{a+3b}{4}\right) }^{2}\lesssim \left\Vert \vartheta
\right\Vert _{L^{2}(a,b)}^{2}+\left\vert \int_{0}^{1}\chi ^{2}\mathcal{P}%
\vartheta \vartheta \right\vert \\
& \lesssim \left\Vert \vartheta \right\Vert _{L^{2}(a,b)}^{2}+\left\vert
\int_{0}^{1}\chi ^{2}\left( \mathcal{P}-\sigma ^{2}\right) \vartheta
\vartheta \right\vert +\sigma ^{2}\left\Vert \vartheta \right\Vert
_{L^{2}(a,b)}^{2} \\
& \lesssim \left\Vert \left( \mathcal{P}-\sigma ^{2}\right) \vartheta
\right\Vert _{L^{2}(0,1)}^{2}+(1+\sigma ^{2})\left\Vert \vartheta
\right\Vert _{L^{2}(a,b)}^{2}
\end{align*}%
by Cauchy-Schwarz. Combining the above estimates ends the proof Proposition
5.1.

\bigskip

Proof of Lemma 5.1 .- Let us consider $\phi (x)=x^{2-\alpha }$ and $%
v=e^{\phi }\vartheta $. Note that for $\alpha \in \lbrack 0,1)$, $v\in D(%
\mathcal{P})$ because $\vartheta \in D(\mathcal{P})$. We set 
\begin{equation*}
P_{\phi }=e^{\phi }\mathcal{P}e^{-\phi }-\sigma ^{2}
\end{equation*}%
with 
\begin{equation*}
\mathcal{S}=-\frac{d}{dx}\left( x^{\alpha }\frac{d}{dx}\right) -(2-\alpha
)^{2}x^{2-\alpha }-\sigma ^{2},\quad \mathcal{A}=2(2-\alpha )x\frac{d}{dx}%
+(2-\alpha )\text{ , }
\end{equation*}%
in order that $P_{\phi }v=e^{\phi }\left( \mathcal{P}-\sigma ^{2}\right)
\vartheta $ and $P_{\phi }v=\mathcal{S}v+\mathcal{A}v$ which gives $%
\left\Vert P_{\phi }v\right\Vert _{L^{2}(0,1)}^{2}=\left\Vert \mathcal{S}%
v\right\Vert _{L^{2}(0,1)}^{2}+\left\Vert \mathcal{A}v\right\Vert
_{L^{2}(0,1)}^{2}+2(\mathcal{S}v,\mathcal{A}v)_{L^{2}(0,1)}$.

\bigskip

Classical computations lead to%
\begin{equation*}
\begin{array}{ll}
(\mathcal{S}v,\mathcal{A}v)_{L^{2}(0,1)} & =(2-\alpha )^{2}\left\Vert
x^{\alpha /2}v^{\prime }\right\Vert _{L^{2}(0,1)}^{2}+(2-\alpha
)^{4}\left\Vert x^{(2-\alpha )/2}v\right\Vert _{L^{2}(0,1)}^{2}-(2-\alpha
)\left\vert v^{\prime }(1)\right\vert ^{2} \\ 
& \text{\quad }+(2-\alpha )\underset{x\rightarrow 0^{+}}{\text{lim}}\left[
x^{1+\alpha }|v^{\prime }(x)|^{2}+x^{\alpha }v^{\prime }(x)v(x)+\left(
2-\alpha \right) ^{2}x^{3-\alpha }\left\vert v(x)\right\vert ^{2}+\sigma
^{2}x\left\vert v(x)\right\vert ^{2}\right] \text{ .}%
\end{array}%
\end{equation*}%
The above limit vanishes from the boundary conditions and the regularity of $%
v$. Therefore, the fact that $0\leq \left\Vert P_{\phi }v\right\Vert
_{L^{2}(0,1)}^{2}-2(\mathcal{S}v,\mathcal{A}v)_{L^{2}(0,1)}$ implies 
\begin{eqnarray*}
2(2-\alpha )^{2}\left\Vert x^{\alpha /2}v^{\prime }\right\Vert
_{L^{2}(0,1)}^{2}+2(2-\alpha )^{4}\left\Vert x^{(2-\alpha )/2}v\right\Vert
_{L^{2}(0,1)}^{2} &\leq &\left\Vert P_{\phi }v\right\Vert
_{L^{2}(0,1)}^{2}+2(2-\alpha )\left\vert v^{\prime }(1)\right\vert ^{2} \\
&=&\left\Vert e^{\phi }\left( \mathcal{P}-\sigma ^{2}\right) \vartheta
\right\Vert _{L^{2}(0,1)}^{2}+2(2-\alpha )\left\vert \vartheta ^{\prime
}(1)\right\vert ^{2}\text{ .}
\end{eqnarray*}%
Since $x^{\alpha /2}\vartheta ^{\prime }=e^{-\phi }\left( x^{\alpha
/2}v^{\prime }-(2-\alpha )x^{\left( 2-\alpha \right) /2}v\right) $, 
\begin{equation*}
\left\Vert x^{\alpha /2}\vartheta ^{\prime }\right\Vert
_{L^{2}(0,1)}^{2}\leq 2\left\Vert x^{\alpha /2}v^{\prime }\right\Vert
_{L^{2}(0,1)}^{2}+2(2-\alpha )^{2}\left\Vert x^{(2-\alpha )/2}v\right\Vert
_{L^{2}(0,1)}^{2}\text{ .}
\end{equation*}%
Combining the two above inequalities, we get, for $\alpha \in \lbrack 0,1)$ 
\begin{eqnarray*}
\left\Vert x^{\alpha /2}\vartheta ^{\prime }\right\Vert _{L^{2}(0,1)}^{2}
&\leq &\frac{1}{(2-\alpha )^{2}}\left( \left\Vert e^{\phi }\left( \mathcal{P}%
-\sigma ^{2}\right) \vartheta \right\Vert _{L^{2}(0,1)}^{2}+2(2-\alpha
)\left\vert \vartheta ^{\prime }(1)\right\vert ^{2}\right) \\
&\lesssim &\left\Vert \left( \mathcal{P}-\sigma ^{2}\right) \vartheta
\right\Vert _{L^{2}(0,1)}^{2}+\left\vert \vartheta ^{\prime }(1)\right\vert
^{2}\text{ .}
\end{eqnarray*}%
It remains to bound $\sigma ^{2}\left\Vert \vartheta \right\Vert
_{L^{2}(0,1)}^{2}$. By Cauchy-Schwarz, 
\begin{eqnarray*}
\sigma ^{2}\left\Vert \vartheta \right\Vert _{L^{2}(0,1)}^{2} &=&\int_{0}^{1}%
\mathcal{P}\vartheta \vartheta -\int_{0}^{1}(\mathcal{P}-\sigma
^{2})\vartheta \vartheta =\left\Vert x^{\alpha /2}\vartheta ^{\prime
}\right\Vert _{L^{2}(0,1)}^{2}-\int_{0}^{1}(\mathcal{P}-\sigma
^{2})\vartheta \vartheta \\
&\leq &\left\Vert x^{\alpha /2}\vartheta ^{\prime }\right\Vert
_{L^{2}(0,1)}^{2}+\left\Vert \vartheta \right\Vert _{L^{2}(0,1)}\left\Vert
\left( \mathcal{P}-\sigma ^{2}\right) \vartheta \right\Vert _{L^{2}(0,1)} \\
&\leq &\left\Vert x^{\alpha /2}\vartheta ^{\prime }\right\Vert
_{L^{2}(0,1)}^{2}+2\left\Vert x\vartheta ^{\prime }\right\Vert
_{L^{2}(0,1)}\left\Vert \left( \mathcal{P}-\sigma ^{2}\right) \vartheta
\right\Vert _{L^{2}(0,1)} \\
&\lesssim &\left\Vert x^{\alpha /2}\vartheta ^{\prime }\right\Vert
_{L^{2}(0,1)}^{2}+\left\Vert \left( \mathcal{P}-\sigma ^{2}\right) \vartheta
\right\Vert _{L^{2}(0,1)}^{2}
\end{eqnarray*}%
where $\left\Vert \vartheta \right\Vert _{L^{2}(0,1)}^{2}\leq 4\left\Vert
x\vartheta ^{\prime }\right\Vert _{L^{2}(0,1)}^{2}$ comes from one
integration by parts. This ends the proof of Lemma 5.1.

\bigskip

Proof of Lemma 5.2 .- Denote $\widetilde{a}=\frac{3a+b}{4}$ and $\widetilde{b%
}=\frac{a+3b}{4}$ in order that $0<a<\widetilde{a}<\widetilde{b}<b<1$. Let
us consider $\phi (x)=e^{\lambda \psi }$, with $\lambda >0$, $\psi \in
C^{\infty }\left( 0,1\right) $, $\psi ^{\prime }\neq 0$ on $[\widetilde{a}%
,1] $ and $\psi ^{\prime }(1)<0$. Let $\chi \in C^{\infty }\left( 0,1\right) 
$ such that $0\leq \chi \leq 1$, $\chi =0$ on $[0,\widetilde{a}]$ and $\chi
=1$ on $[\widetilde{b},1]$, and let $v=e^{\tau \phi }\chi \vartheta $ with $%
\tau >0$. We set 
\begin{equation*}
P_{\phi }=e^{\tau \phi }\mathcal{P}e^{-\tau \phi }-\sigma ^{2}
\end{equation*}%
with 
\begin{equation*}
\mathcal{S}=-\frac{d}{dx}\left( x^{\alpha }\frac{d}{dx}\right) -\tau
^{2}x^{\alpha }\left\vert \phi ^{\prime }\right\vert ^{2}-\sigma ^{2},\quad 
\mathcal{A}=2\tau x^{\alpha }\phi ^{\prime }\frac{d}{dx}+\tau \left(
x^{\alpha }\phi ^{\prime }\right) ^{\prime }\text{ ,}
\end{equation*}%
in order that $P_{\phi }v=e^{\tau \phi }\left( \mathcal{P}-\sigma
^{2}\right) \left( \chi \vartheta \right) $ and $P_{\phi }v=\mathcal{S}v+%
\mathcal{A}v$ which gives $\left\Vert P_{\phi }v\right\Vert
_{L^{2}(0,1)}^{2}=\left\Vert \mathcal{S}v\right\Vert
_{L^{2}(0,1)}^{2}+\left\Vert \mathcal{A}v\right\Vert _{L^{2}(0,1)}^{2}+2(%
\mathcal{S}v,\mathcal{A}v)_{L^{2}(0,1)}$.

\bigskip

Classical computations lead to 
\begin{multline*}
(\mathcal{S}v,\mathcal{A}v)_{L^{2}(0,1)}=2\tau \int_{0}^{1}x^{2\alpha }\phi
^{\prime \prime }\left\vert v^{\prime }\right\vert ^{2}+\tau \alpha
\int_{0}^{1}x^{2\alpha -1}\phi ^{\prime }\left\vert v^{\prime }\right\vert
^{2} \\
-\frac{\tau }{2}\int_{0}^{1}\left( \mathcal{P}^{2}\phi \right) \left\vert
v\right\vert ^{2}+2\tau ^{3}\int_{0}^{1}x^{2\alpha }\phi ^{\prime \prime
}\left( \phi ^{\prime }\right) ^{2}\left\vert v\right\vert ^{2}+\alpha \tau
^{3}\int_{0}^{1}x^{2\alpha -1}(\phi ^{\prime })^{3}\left\vert v\right\vert
^{2}-\tau \phi ^{\prime }\left( 1\right) \left\vert v^{\prime }\left(
1\right) \right\vert ^{2}\text{ .}
\end{multline*}%
But, using $\phi =e^{\lambda \psi }$ with $\psi $ having a non-vanishing
gradient, there exist five constants $C_{0},C_{1},C_{2},C_{3},C_{4}>0$
independent on $\alpha \in \lbrack 0,1)$ such that 
\begin{multline*}
(\mathcal{S}v,\mathcal{A}v)_{L^{2}(0,1)}\geq \tau \lambda
^{2}C_{0}\int_{0}^{1}\phi \left\vert v^{\prime }\right\vert ^{2}+\tau
^{3}\lambda ^{4}C_{1}\int_{0}^{1}\phi ^{3}\left\vert v\right\vert ^{2}-\tau
\lambda C_{2}\int_{0}^{1}\phi \left\vert v^{\prime }\right\vert ^{2} \\
-\tau \lambda ^{4}C_{3}\int_{0}^{1}\phi \left\vert v\right\vert ^{2}-\tau
^{3}\lambda ^{3}C_{4}\int_{0}^{1}\phi ^{3}\left\vert v\right\vert ^{2}+\tau
\left\vert \phi ^{\prime }\left( 1\right) \right\vert \left\vert v^{\prime
}\left( 1\right) \right\vert ^{2}\text{ .}
\end{multline*}%
Therefore, the fact that $0\leq \left\Vert P_{\phi }v\right\Vert
_{L^{2}(0,1)}^{2}-2(\mathcal{S}v,\mathcal{A}v)_{L^{2}(0,1)}$ implies by
taking $\lambda >0$ sufficiently large, and $\tau >0$ sufficiently large the
following inequality%
\begin{equation*}
\left\Vert v\right\Vert _{H^{1}(0,1)}^{2}+\left\vert v^{\prime }\left(
1\right) \right\vert ^{2}\lesssim \left\Vert P_{\phi }v\right\Vert
_{L^{2}(0,1)}^{2}\text{ .}
\end{equation*}%
Taking the weights off the integrals and using commutators, we have%
\begin{equation*}
|\vartheta ^{\prime }\left( 1\right) |^{2}\lesssim \left\Vert P_{\phi
}\vartheta \right\Vert _{L^{2}(0,1)}^{2}+\left\Vert \vartheta \right\Vert
_{H^{1}(\widetilde{a},\widetilde{b})}^{2}\text{ .}
\end{equation*}%
This ends the proof of Lemma 5.2.

\bigskip

\section{Observability estimate for the degenerate heat equation (proof of
Theorem 1.2)}

\bigskip

In this section, we prove that the refine observability from measurable set
of Theorem 1.2 is a corollary of the spectral Lebeau-Robbiano inequality of
Theorem 1.1.

\bigskip

Let $\widetilde{\omega }\Subset \omega $ and $\chi \in C_{0}^{\infty }\left(
\omega \right) $ be such that $0\leq \chi \leq 1$ and $\chi =1$ in $%
\widetilde{\omega }$.

\bigskip

We start with Theorem 3.1 of \cite[page 1142]{BP} stating that $\left(
i\right) $ implies $\left( ii\right) $ where:

$\left( i\right) $ $\exists C_{1}>0$, $\forall \left\{ a_{j}\right\} \in 
\mathbb{R}$, $\forall \Lambda >0$ 
\begin{equation*}
\sum_{\lambda _{j}\leq \Lambda }|a_{j}|^{2}\leq e^{C_{1}\left( 1+\sqrt{%
\Lambda }\right) }\int_{\widetilde{\omega }}\left\vert \sum_{\lambda
_{j}\leq \Lambda }a_{j}\Phi _{j}\right\vert ^{2}\text{ ;}
\end{equation*}

$\left( ii\right) $ $\forall t>0$, $\forall \varepsilon \in \left(
0,2\right) $, $\forall y_{0}\in L^{2}\left( 0,1\right) $%
\begin{equation*}
\left\Vert e^{-t{\mathcal{P}}}y_{0}\right\Vert _{L^{2}\left( 0,1\right)
}\leq 4e^{\frac{C_{1}}{2}}e^{\frac{C_{1}^{2}}{2\varepsilon t}}\left\Vert
e^{-t{\mathcal{P}}}y_{0}\right\Vert _{L^{2}\left( \widetilde{\omega }\right)
}^{1-\varepsilon /2}\left\Vert y_{0}\right\Vert _{L^{2}\left( 0,1\right)
}^{\varepsilon /2}\text{ .}
\end{equation*}

Therefore, by Theorem 1.1 we know that $\left( ii\right) $ holds with $%
C_{1}=C\frac{1}{\left( 2-\alpha \right) ^{2}}>1$.

\bigskip

By Nash inequality and regularizing effect, we get for some constants $c>1$
and $\theta \in \left( 0,1\right) $ independent on $\left( y_{0},t\right) $
and $\alpha \in \left[ 0,2\right) $ 
\begin{equation*}
\left\Vert e^{-t{\mathcal{P}}}y_{0}\right\Vert _{L^{2}\left( \widetilde{%
\omega }\right) }\leq c\left( 1+\frac{1}{\sqrt{t}}\right) ^{\theta
}\left\Vert e^{-t{\mathcal{P}}}y_{0}\right\Vert _{L^{1}\left( \omega \right)
}^{1-\theta }\left\Vert y_{0}\right\Vert _{L^{2}\left( 0,1\right) }^{\theta }%
\text{ .}
\end{equation*}

\bigskip

Therefore, since $1+\frac{1}{\sqrt{t}}\leq 4e^{\frac{C_{1}^{2}}{2\varepsilon
t}}$, with $C_{2}=16ce^{\frac{C_{1}}{2}}e^{\frac{C_{1}^{2}}{\varepsilon t}%
}\geq 4e^{\frac{C_{1}}{2}}e^{\frac{C_{1}^{2}}{2\varepsilon t}}\left( c\left(
1+\frac{1}{\sqrt{t}}\right) ^{\theta }\right) ^{1-\varepsilon /2}$ and $\mu
=1-\left( 1-\theta \right) \left( 1-\frac{\varepsilon }{2}\right) $, we
obtain 
\begin{equation*}
\left\Vert e^{-t{\mathcal{P}}}y_{0}\right\Vert _{L^{2}\left( 0,1\right)
}\leq C_{2}\left\Vert e^{-t{\mathcal{P}}}y_{0}\right\Vert _{L^{1}\left(
\omega \right) }^{1-\mu }\left\Vert y_{0}\right\Vert _{L^{2}\left(
0,1\right) }^{\mu }\text{ ,}
\end{equation*}%
which implies by Young inequality%
\begin{eqnarray*}
\left\Vert e^{-t{\mathcal{P}}}y_{0}\right\Vert _{L^{2}\left( 0,1\right) }
&\leq &s\left\Vert y_{0}\right\Vert _{L^{2}\left( 0,1\right) }+\frac{1}{s^{%
\frac{\mu }{1-\mu }}}C_{2}^{\frac{1}{1-\mu }}\left\Vert e^{-t{\mathcal{P}}%
}y_{0}\right\Vert _{L^{1}\left( \omega \right) } \\
&\leq &s\left\Vert y_{0}\right\Vert _{L^{2}\left( 0,1\right) }+\frac{1}{s^{%
\frac{\mu }{1-\mu }}}\left( 16ce^{\frac{C_{1}}{2}}\right) ^{\frac{1}{1-\mu }%
}e^{\frac{C_{1}^{2}}{\varepsilon t}\frac{1}{1-\mu }}\left\Vert e^{-t{%
\mathcal{P}}}y_{0}\right\Vert _{L^{1}\left( \omega \right) }\text{ .}
\end{eqnarray*}

\bigskip

Reproducing the proof of Theorem 1.1 of \cite[page 684]{PW}, we have for our
system that $\left( iii\right) $ implies $\left( iv\right) $ where:

$\left( iii\right) $ $\exists K_{1},K_{2},\ell >0$, $\forall s>0$%
\begin{equation*}
\left\Vert e^{-t{\mathcal{P}}}y_{0}\right\Vert _{L^{2}\left( 0,1\right)
}\leq s\left\Vert y_{0}\right\Vert _{L^{2}\left( 0,1\right) }+\frac{1}{%
s^{\ell }}K_{1}e^{\frac{K_{2}}{t}}\left\Vert e^{-t{\mathcal{P}}%
}y_{0}\right\Vert _{L^{1}\left( \omega \right) }\text{ ;}
\end{equation*}

$\left( iv\right) $ $\forall y_{0}\in L^{2}\left( 0,1\right) $%
\begin{equation*}
\left\Vert e^{-T{\mathcal{P}}}y_{0}\right\Vert _{L^{2}\left( 0,1\right)
}\leq K_{3}\int_{\omega \times E}\left\vert e^{-t{\mathcal{P}}%
}y_{0}\right\vert \text{ ,}
\end{equation*}%
with 
\begin{equation*}
K_{3}=c\frac{K_{1}}{K_{2}}e^{cK_{2}}\text{ when }E\subset \left( 0,T\right) 
\text{ is a measurable set of positive measure;}
\end{equation*}%
\begin{equation*}
K_{3}=\kappa \frac{K_{1}}{K_{2}}e^{\kappa \frac{K_{2}}{T}}\text{ when }%
E=\left( 0,T\right) \text{ for some }\kappa >0\text{ independent on }T\text{
.}
\end{equation*}

\bigskip

Therefore, with $K_{1}=\left( 16ce^{\frac{C_{1}}{2}}\right) ^{\frac{1}{%
\left( 1-\theta \right) \left( 1-\frac{\varepsilon }{2}\right) }}$ and $%
K_{2}=\frac{1}{\varepsilon \left( 1-\theta \right) \left( 1-\frac{%
\varepsilon }{2}\right) }C_{1}^{2}$ we have $K_{3}\leq Ce^{C\frac{1}{%
(2-\alpha )^{4}}}$. This completes the proof of Theorem 1.2.

\bigskip

\bigskip

\bigskip

\end{document}